%% file: Floer_vs_categorical.tex
\documentclass[11pt]{amsart}
\usepackage{amsmath}
\usepackage{amsfonts}
\usepackage{amssymb}
\usepackage{graphicx}
\usepackage{hyperref}
\usepackage{mathrsfs}
\usepackage{pb-diagram}
\usepackage{epstopdf}
\usepackage{amsmath,amsfonts,amsthm,enumerate,amscd,latexsym,curves}
\usepackage{bbm}
\usepackage[mathscr]{eucal}
\usepackage{epsfig,epsf,epic}
\usepackage{caption}
\usepackage{subcaption}
\usepackage{multirow}
\usepackage{color}
\usepackage[all]{xy}
\usepackage{fourier}
\usepackage{tikz-cd}
\usepackage{calc}
\usepackage{float}
\usepackage{enumitem}
\usepackage{comment}
\usepackage{faktor}

\newcommand{\HF}{\mathrm{HF}}
\newcommand{\hcat}{h_{\mathrm{cat}}}
\newcommand{\hfl}{h_{\mathrm{Fl}}}

\newcommand{\tr}{\mathrm{tr}\left(M_f\right)}
\newcommand{\arccosh}{\mathrm{arccosh}}

\newtheorem{thm}{Theorem}[section]
\newtheorem{lem}[thm]{Lemma}

\newtheorem{cor}[thm]{Corollary}
\newtheorem{prop}[thm]{Proposition}
\newtheorem{definition}[thm]{Definition}

\newtheorem{remark}[thm]{Remark}

\numberwithin{equation}{section}

\begin{document}

\title[A comparison of entropies and counting conjugacy classes]{Categorical and Floer-theoretic entropies: comparison and counting conjugacy classes in symplectic mapping class groups}

\author{Sangjin Lee}
\address[Sangjin Lee]{Korea Institute for Advanced Study}
\email{sangjinlee@kias.re.kr}

\maketitle
%\tableofcontents

%=================================================================================================================%
\begin{abstract}
	Let $W$ be a Weinstein manifold and let $f: W \to W$ be a compactly supported symplectic automorphism. 
	Using Floer theories, one can define two different entropies of $f$, called {\em Floer-theoretic entropy} and {\em categorical entropy}. 
	The current article compares these two entropies and show that the former bounds the latter from below.
	 
	As an application of the comparison, using that these entropies are invariants up to the conjugate relation, we count conjugacy classes in symplectic mapping classes. 
	Especially, we show an exponential abundance of conjugacy classes that are smoothly isotopic to the identity, when $W$ is a Weinstein manifold of dimension $4n \geq 6$, satisfying a local condition.
\end{abstract}

\section{Introduction}
\label{section introduction}
Classical mechanics studies motions of moving particles, thus it motivates interesting dynamics problems, such as the three-body problem.
Since symplectic topology has been originated from classical mechanics, symplectic topology could help us to understand these dynamics problems. 
Especially, one can formulate a large class of classical mechanical systems as {\em Hamiltonian dynamics} defined on symplectic manifolds.  
It motivates to study Hamiltonian dynamical systems in symplectic topology. 

To study Hamiltonian dynamical systems, Poincar\'e suggested an idea to use {\em Poincar\'e's return map}. 
More specifically, Poincar\'e's idea is to study Hamiltonian dynamics by recording its intersections with a lower-dimensional symplectic section.
On the section, Hamiltonian dynamics becomes a symplectic automorphism of the section, which is not necessarily to be a Hamiltonian dynamics again. 
We note that every symplectic automorphism can occur as a Poincar\'e return map of a Hamiltonian dynamics. 
It motivates to study symplectic automorphisms and their dynamics, or simply, {\em symplectic dynamical systems}.

In the article, we study symplectic dynamical systems. 
By studying their dynamical properties, we can find exponentially abundant symplectic automorphisms that are smoothly isotopic to the identity, on Liouville manifolds under a condition.

To specify our theme and to provide more details, we discuss the following two practical questions: 
First, what aspects/properties of symplectic dynamical systems do we want to study? Second, in order to study them, what tools can we use? 

The answer for the first question is {\em dynamical complexities} of symplectic dynamics. 
Dynamical complexity is one of the main topics in Dynamics. 
To record the complexity, dynamicists defined the notion of {\em entropy} of a dynamical system: it is a non-negative real number assigned to a dynamical system such that, if a dynamical system has positive entropy, it has topologically chaotic behavior.
For more details on this basic notion in Dynamics, one can consult textbooks such as \cite{Brin-Stuck15}.  

Since we are interested in {\em symplectic} dynamical systems, we would like to measure their complexities via symplectic topological tools.
Specifically, we would like to employ various {\em Floer theories} in symplectic topology. 
Let us note that historically, Floer theories have been originally utilized to study symplectic dynamics: 
Floer \cite{Floer87} introduced the first Floer homology, {\em fixed-point Floer homology}, in order to study periodic orbits of Hamiltonian dynamics. 
Also, in \cite{Floer88}, he defined the {\em Lagrangian Floer homology} and as an application of it, he could study the Hamiltonian fixed-point problem.
For more explanations, we refer the reader to \cite{Abbondandolo-Schlenk19}.

Floer theories have been used to record dynamical complexities of symplectic automorphisms. 
For example, Fel'shtyn \cite{Felshtyn04} defined an {\em asymptotic invariant} of a symplectic automorphism $f$ by measuring the exponential growth of $\dim \HF_{fixed}\left(f^n\right)$, where $\HF_{fixed}\left(f^n\right)$ is the fixed-point Floer homology of $f^n$. 
Similarly, Smith \cite{Smith12} introduced the following logarithmic definition
\begin{gather}
	\label{eqn Floer theoretic entropy}
	\hfl(f) := \limsup_{n \to \infty} \frac{1}{n} \log \dim \HF_{fixed} \left(f^n\right),
\end{gather}
and named it {\em Floer-theoretic entropy}. 

We note that counting periodic points of a dynamical system $f$, i.e., fixed-points of $f^n$ for some $n \in \mathbb{N}$, is a natural way to measure the complexity of $f$ in Dynamics.
For example, see \cite{Artin-Mazur65, Bowen70}.
Since, by definition of fixed-point Floer homology, $\hfl(f)$ is a lower bound of the exponential growth of a number of periodic points of $f$, $\hfl(f)$ could measure some dynamical complexity of $f$. 
It is summarized in \cite{Smith12}: Smith described $\hfl(f)$ as a ``robust version of the periodic entropy.'' 
(We also note that \cite{Smith12} used positive Floer-theoretic entropy to study the faithfulness of mapping class group actions.
It motivates our application explained below.)

Another example of using Floer theory to measure complexities of symplectic dynamical systems is the {\em categorical entropy}:
provided a symplectic dynamical system $f:W \to W$, Lagrangian Floer theory defines the {\em Fukaya category} of $W$.
Also, the Floer theory explains how $f$ induces a functor $\phi$ on the Fukaya category of $W$. 
Then, the notion of categorical entropy defined in \cite{Dimitrov-Haiden-Katzarkov-Kontsevich} measures the complexity of the categorical dynamical system $\phi$, which can be seen as an invariant of the original $f$. 
Let this invariant be called {\em categorical entropy of $f$} and let us denote it by $\hcat(f)$. 

Let us explain briefly why $\hcat(f)$ measures the (topological) dynamical complexity of $f$.
In the article, we work in the setting of Weinstein manifolds.
It means that there exits a finite number of Lagrangian submanifolds that {\em generate} every Lagrangian. 
Roughly speaking, it means that every Lagrangian can be constructed by taking surgeries of copies of generating Lagrangians. 
The categorical entropy $\hcat(f)$ cares the number of copies of generators that we need to generate $f^n(L)$ for Lagrangians $L$.
More specifically, (but still imprecisely,)
\begin{gather}
	\label{eqn categorical entropy}
	\hcat(f) = \max_L \lim_{n \to \infty} \frac{1}{n} \log \left(\text{the minimum number of generators we need to generate  } f^n(L) \right).
\end{gather}  
We note that if $L$ is more ``complicated'' Lagrangian, then we need more copies of generating Lagrangians in order to generate $L$. 
Thus, one can speak that $\hcat(f)$ measures the complexity of the action of $f$ on Lagrangians. 

Interestingly, it is known in the literature that $\hfl(f)$ and $\hcat(f)$ help us to measure the dynamical complexity of $f$, in the sense of {\em Dynamics}, not symplectic topology or Floer theories. 
To be more precise, let $h(f)$ denote the topological entropy that Dynamics studies. 
Then, the following inequalities are known in \cite{Bae-Lee25, Kim-Kwon26} (for the case of Weinstein manifolds $W$ and compactly supported $f$):
\[\hcat(f), \hfl(f) \leq h(f).\]

The known inequalities suggest a possible application in the theory of dynamical systems. 
As mentioned above, the study of $h(f)$ occupies a central role in the theory.
Especially, by showing positivity of $h(f)$, one can say that $f$ possesses topologically chaotic behavior. 
Unfortunately, direct computations of $h(f)$ are often unachievable.
However, the inequalities can detect positivity of $h(f)$ by utilizing Floer theories and symplectic topology, which topological methods, such as homology, cannot detect.  
See \cite[Section 5]{Bae-Choa-Jeong-Karabas-Lee25} for example. 

Since we have two different entropies of a symplectic dynamical system, it would be natural request to compare them.
More specific question would be ``what is the relation between two complexities of a symplectic dynamical system, detected from fixed-point and Lagrangian Floer theories?''
Theorem \ref{thm intro} answers the question. 
\begin{thm}
	\label{thm intro}
	Let $W$ be a Weinstein manifold of finite type, and let $f:W \to W$ be a compactly supported, exact symplectic automorphism of $W$. 
	Then,
	\[\hfl(f) \leq \hcat(f).\] 
\end{thm}

\begin{proof}[Proof sketch]
	We briefly explain our idea for Theorem \ref{thm intro} in the introduction section. 
	The full proof will be stated in Section \ref{section proof of the main theorem}. 
	
	The difficulty of proving Theorem \ref{thm intro} arises from the fact that the notion of $\hfl(f)$ and $\hcat(f)$ utilize {\em different} Floer theories. 
	Our idea is to understand them in a single Floer theory. 
	The single Floer theory is {\em a version} (will be explained soon) of Lagrangian Floer homology of the product space $W \times W$. 
	
	The comparison between Fixed-point Floer homology and Lagrangian Floer homology of the product is conventional: 
	Floer \cite{Floer88} related $\HF_{fixed}(f)$ and $\HF_{\left(W \times W\right)}\left(\Delta, \Gamma(f)\right)$, where the latter means Lagrangian Floer homology in $W \times W$ between the diagonal $\Delta$ and the graph $\Gamma(f)$ defined as follows: 
	\[\Delta := \left\{(x,x) \in W \times W : x \in W\right\}, \quad \Gamma(f) := \left\{(x,f(x)) \in W \times W : x \in W\right\}.\]
	
	We note that in \cite{Kim-Kwon26}, the authors studied the conventional comparison in our setting, i.e., the setting of Weinstein manifold $W$ and compactly supported $f$. 
	More specifically, \cite[Theorem 2.6]{Kim-Kwon26} showed that the fixed-point Floer homology of $f$ is equivalent to a {\em partially wrapped} Lagrangian Floer homology between $\Delta$ and $\Gamma(f)$, i.e.,
	\[\HF_{fixed}(f) \simeq \HF_{\left(W \times W\right)}\left(\Delta, \Gamma(f)\right),\]
	where the later means the partially wrapped Lagrangian Floer homology by abuse of notation.  
	
	We note that 
	\[\Gamma(f^n) = F^n\left(\Delta\right), \text{  where  } F: W \times W \to W \times W, \quad (x,y) \mapsto \left(x,f(y)\right).\]
	Thus, the above arguments can be summarized in the following way: 
	$\hfl(f)$ measures the complexity of $F$ with respect to a {\em specific} Lagrangian $\Delta$. 
	
	By Equation \eqref{eqn categorical entropy}, $\hcat(F)$ measures complexity of $F$ with respect to {\em all} Lagrangians of $W \times W$. 
	Thus, it is natural to expect that $\hfl(f) \leq \hcat(F)$. 
	And, by the K\"unneth formula from \cite{Ganatra-Pardon-Shende24}, we could expect that $\hcat(F) \leq \hcat(f)$. 
	Section \ref{section proof of the main theorem} proves that these expectations hold, and it completes the proof.
\end{proof}

\begin{remark}
	\label{rmk barcode entropy}
	We would like to mention that there exists one more entropy type invariant, defined for symplectic dynamical systems, by employing Floer theory. 
	The extra one is called {\em barcode entropy}.
	See \cite{Cineli-Ginzburg-Gurel24} for the definition.
	In \cite{Bae-Lee25}, under the setting of Theorem \ref{thm intro}, it is known that $\hcat(f)$ is smaller than or equal to the barcode entropy of $f$, and this comparison has been used to show that $\hcat(f)$ is smaller than or equal to the topological entropy of $f$. 
	In summary, we have 
	\[\hfl(f) \leq \hcat(f) \leq \text{  the barcode entropy of  } f \leq h(f).\]
\end{remark}

A natural question follows Theorem \ref{thm intro}: 
how sharp is the inequality?
In the current article, we do not fully answer the question, but we give two examples helping us to understand the question. 

The first example, explained in Section \ref{subsection the torus example}, shows that the strict inequality can happens, but in a different setting from Theorem \ref{thm intro}.
More specifically, we give a symplectic automorphism $f$ on a {\em closed} symplectic manifold (not a Weinstein manifold) $\mathbb{T} \times \mathbb{T}$, where $\mathbb{T}$ is the $2$-dimensional torus, and show that 
\[0 = \hfl(f) < \hcat(f).\] 
We note again that, it happens in a different setting from that of Theorem \ref{thm intro} and we do not know if the strict inequality can happen in the original setting.

The second example, given in Section \ref{subsection symplectic automorphisms of Penner type}, is the symplectic automorphisms of {\em Penner type} (see Definition \ref{definition Penner type}). 
Their categorical entropy have been computed in \cite{Bae-Choa-Jeong-Karabas-Lee25} and \cite{Kim-Kwon26} computed a lower-bound for their Floer-theoretic entropy.
Interestingly, one can see that the lower-bound from \cite{Kim-Kwon26} coincides with the categorical entropy from \cite{Bae-Choa-Jeong-Karabas-Lee25}.
Thus, $\hfl(f)=\hcat(f)$ if $f$ is of Penner type, and it shows that the inequality could be sharp at least for some symplectic dynamics.  
The example can be seen as a computational application of Theorem \ref{thm intro}, which computes the exact Floer-theoretic entropy.

Employing Theorem \ref{thm intro}, we study {\em symplectic mapping class groups}. 
To state our result more precisely, we define notation first. 
For a Weinstein manifold $W$, let 
\[\textit{Symp}(W) := \left\{f: W \to W : f \text{  is a compactly supported, exact symplectic automorphism}\right\}.\]
Then, $\textit{Symp}(W)$ admits a natural $C^\infty$-topology, and the {\em symplectic mapping class group of $W$} is defined to be $\pi_0\left(\textit{Symp}(W)\right)$. 

A natural subgroup of $\pi_0\left(\textit{Symp}(W)\right)$ is the {\em smoothly trivial symplectic mapping class group}, defined by 
\[K(W):= \ker \left[\pi_0\left(\textit{Symp}(W)\right) \xrightarrow{\text{forgetful map}} \pi_0\left(\textit{Diff}^+(W)\right)\right],\]
where $\textit{Diff}^+(W)$ is the group of orientation preserving diffeomorphisms of $W$.
When $W$ is two-dimensional, Moser's theorem implies that $\pi_0\left(\textit{Symp}(W)\right)$ and $\pi_0\left(\textit{Diff}^+(W)\right)$ are identical by the forgetful map, and hence $K(W)$ is trivial. 
In higher-dimensional, however, $K(W)$ can be nontrivial, as proven in many references including \cite{Seidel99, Khovanov-Seidel02, Keating14, Keating22}.
In other words, there exist symplectic automorphisms that are smoothly isotopic to the identity, but not symplectically so.
It elucidates the difference between smooth topology and symplectic topology.
It motivates to study $K(W)$. 
 
As an application of Theorem \ref{thm intro}, we measure the size of $K(W)$ with respect to the Floer-theoretic entropy, after passing to conjugacy classes.
For a formal statement of our result, we define the conjugate equivalence relation for $f, g \in \pi_0 \left(\textit{Symp}(W)\right)$ by
\[f \sim g \iff \text{  there exists  } h \in \pi_0 \left(\textit{Symp}(W)\right) \text{  such that  } f = h g h^{-1}.\]
We give a name for the quotient by the conjugate relation:
\[C(W):= \faktor{\pi_0 \left(\textit{Symp}(W)\right)}{\sim} = \left\{[f] : f \in \pi_0 \left(\textit{Symp}(W)\right)\right\}.\]
Moreover, since if $f \in K(W)$, then for all $g \sim f$, $g \in K(W)$, i.e., $K(W)$ is normal, one can define 
\[C_{tri}(W):= \faktor{K(W)}{\sim} = \left\{[f] \in \left(\faktor{\pi_0\left(\textit{Symp}(W)\right)}{\sim}\right) \Big\vert f \in K(W)\right\}.\]
Moreover, we point out that the Floer-theoretic entropy is an invariant of conjugacy classes, see \cite{Smith12}.
Using them, our result can be stated formally as follows:
\begin{thm}
	\label{thm application intro}
	Let $W$ be a Liouville manifold of dimension $\geq 6$ such that there exist two exact Lagrangian spheres intersecting transversally at one point. 
	\begin{enumerate}
		\item The size of the set 
		\[C\left(W;H\right) := \left\{ [f] \in C(W) : \hfl(f) <H\right\}\]
		is bounded from below by an exponential function of $H$. 
		In other words, there exist exponentially abundant conjugacy classes in $\pi_0\left(\textit{Symp}(W)\right)$.
		\item If $\dim W = 4n \geq 6$, the size of the set
		\[C_{tri}\left(W;H\right) := \left\{ [f] \in C_{tri}(W) : \hfl(f) < H\right\},\]
		is bounded from below by an exponential function of $H$. 
		In other words, there exist exponentially abundant conjugacy classes in $K(W)$.
	\end{enumerate} 
\end{thm}

Before sketching the proof of Theorem \ref{thm application intro}, we explain a motivation for counting {\em conjugacy classes} in $\pi_0\left(\textit{Symp}(W)\right)$, a meaning of Theorem \ref{thm application intro}, as well as where Theorem \ref{thm application intro} is located in the literature. 

A motivation arises from surface theory.
Let us recall that for a surface, its symplectic mapping class group and the usual mapping class group are identical. 
And, the conjugacy classes of the mapping classes have strong relation to the {\em Moduli space} of the surface. 

To be more precise, let $\Sigma_g$ denote the closed surface of genus $g$. 
Then, the {\em Moduli space} of $\Sigma_g$ is defined by the Teichm\"uller space quotient by the mapping class group of it, where Teichm\"uller space is the 
Teichm\"uller space of $\Sigma_g$, the space of parameterized complex structure on $\Sigma_g$.
We note that the Moduli space admits a metric, called {\em Teichm\"uller metric}.
In the setting, \cite{Eskin-Mirzakhani11} counted the number of {\em certain} conjugacy classes of the mapping classes, having bounded {\em stretching factor} $\lambda$, 
and the result was 
\[\left\vert \{\text{pseudo-Anosov conjugacy classes  } [f] | f \in \pi_0 \left(\textit{Diff}^+\left(\Sigma_g\right), \text{  the stretching factor of  } f \leq e^H\right\}\right\vert \sim \frac{e^{(6g-g)H}}{(6g-g)H}\]
The number coincides with the number of {\em Teichmüller geodesics} of the Moduli space. 

To discus its relation to symplectic topology, we note that there established connections between Teichm\"uller theory and the space of stability conditions on Fukaya categories, by  Gaiotto--Moore--Neitzke \cite{Gaiotto-Moore-Neitzke13}, Bridgeland--Smith \cite{Bridgeland-Smith15}, Haiden-Katzarkov-Kontsevich \cite{Haiden-Katzarkov-Kontsevich17}, Dimitrov--Haiden--Katzarkov--Kontsevich \cite{Dimitrov-Haiden-Katzarkov-Kontsevich}, Fan--Filip--Haiden--Katzarkov--Liu \cite{Fan-Filip-Haiden-Katzarkov-Liu21} and many others.  
Moreover, Fel'shtyn \cite{Felshtyn12} shows the equivalence between Floer-theoretic entropy and the stretching factor for pseudo-Anosov surface mapping classes. 
These motivate the counting of conjugacy classes of symplectic mapping classes, having bounded Floer-theoretic entropy.
Theorem \ref{thm application intro} (1) does the counting and shows the {\em at-least-exponential} abundance of conjugacy classes.

\begin{remark}
	\label{rmk pseudo-Anosov}
	Compared to \cite{Eskin-Mirzakhani11}, Theorem \ref{thm application intro} omitted the pseudo-Anosov condition.
	There are two possible excuses:
	first, it is known that pseudo-Anosov surface mapping classes are {\em generic} in surface mapping classes. 
	Thus, the counting in \cite{Eskin-Mirzakhani11} could be seen as a counting of {\em generic} conjugacy classes and one may expect that for all conjugacy classes, the same exponential growth rate holds. 
	Another reason why we omit the pseudo-Anosov condition is that there is no higher-dimensional analogy of pseudo-Anosov, from the viewpoint of symplectic topology. 
	However, there exists a {\em categorical} generalization of pseudo-Anosov, see \cite{Fan-Filip-Haiden-Katzarkov-Liu21, Bae-Lee23}.
	And, actually, in the proof of Theorem \ref{thm application intro}, we show the exponential abundance of conjugacy classes of symplectic mapping classes that induce {\em pseudo-Anosov} auto-equivalences on the corresponding category.
\end{remark}

Another motivation for counting conjugacy classes arises from the study of $K(W)$. 
Since a single conjugacy class $[f] \in C(W)$ can have infinitely many distinct mapping classes; such examples can be found in \cite{Keating22}. 
Thus, by showing the exponential lower-bound in Theorem \ref{thm application intro}, one can show that the richness of $K(W)$ can survive even after passing to the quotient by the conjugate relation. 

We point out that Theorem \ref{thm application intro} assumes a local condition.
An interesting point of Theorem \ref{thm application intro} is that the local condition forces the symplectic mapping class group of the {\em ambient} $W$, i.e., $\pi_0 \left(\textit{Symp}(W)\right)$, to have exponentially abundant conjugacy classes, without studying the global structure of $\pi_0 \left(\textit{Symp}(W)\right)$. 
It makes Theorem \ref{thm application intro} easy to be applied.

Moreover, if we can replace the local condition with another, Theorem \ref{thm application intro} could be more applicable:
for example, one can consider two Lagrangian spheres transversally intersecting at {\em two or more} points, instead of them intersecting at one point. 
Remark \ref{rmk local condition} discusses it briefly.
%A_2 plumbing 대신 Kronecker quiver plumbing 등을 넣을 수 있다는 remark필요 - grading에 대한 condition도 들어가야 함! 

In order to clarify where Theorem \ref{thm application intro} is located in the literature, we list some related results, especially related to the counting of conjugacy classes of smoothly trivial mapping classes: 
\begin{itemize}
	\item Keating \cite{Keating14} shows that under a local condition on two Lagrangian spheres, their Dehn twists generate a free subgroup of the symplectic mapping class group.
	It also shows that under a dimension condition, $K(W)$ admits a free subgroup. 
	\item Keating \cite{Keating22} shows that Right-Angled Artin Group $A$ can be embedded into $\pi_0 \left(\textit{Symp}(W)\right)$ (or $K(W)$) for some $W$, see \cite[Corollary 3.11 and Remark 3.12]{Keating22}. 
	It implies that $\pi_0 \left(\textit{Symp}(W)\right)$ can possess finitely presented subgroups $H$ with unsolvable conjugacy problems. 
	(The unsolvable conjugacy problem concerns the subgroup $H$, not $\pi_0 \left(\textit{Symp}(W)\right)$.)
	\item In \cite{Keating-Smith24, Keating-Smith-Wemyss26}, authors show that symplectic mapping class groups of certain Weinstein manifolds can {\em split}.
	Using these splittings, one can study conjugacy classes of $\pi_0 \left(\textit{Symp}(W)\right)$.  
	We note that \cite{Keating-Smith24, Keating-Smith-Wemyss26} employed Fukaya categories of Weinstein manifolds.
	It means that these results do not assume local conditions, differently from Theorem \ref{thm application intro}.
	\item Recently, Kim--Kwon \cite{Kim-Kwon26} shows that $K(W)$ has infinitely many mapping classes, if $W$ satisfies a local condition. 
	Also, using their method, one can show that $C(W)$ also has infinitely many conjugacy classes. 
	However, \cite{Kim-Kwon26} does not show the exponential abundance. 
\end{itemize}

Finally, we sketch the proof of Theorem \ref{thm application intro}.
\begin{proof}[Proof sketch of Theorem \ref{thm application intro}]
	The local condition, i.e., the existence of two Lagrangian spheres intersecting transversally at one point, implies that the $A_2$-Milnor fiber is a Liouville-subdomain of $W$, see Lemma \ref{lem subdomain}. 
	Let $P$ denote the subdomain inside $W$.
	Then, if $f:P \to P$ is an exact symplectic automorphism supported in the interior of $P$, $f$ can extend to the whole Weinstein manifold $W$ naturally, by setting as the identity outside of $P$. 
	Let $\tilde{f}$ denote the extended symplectic automorphism on $W$. 
	
	The local-to-global property of Floer-theoretic entropy, \cite[Theorem 3.2]{Kim-Kwon26}, shows that $\hfl(f) = \hfl(\tilde{f})$. 
	We note that the former is computed in $P$ and the latter is computed in $W$. 
	Thus, in order to prove the first part, it is enough, for all $H \in \mathbb{R}_{>0}$, to construct a finite collection of symplectic automorphisms, 
	\[\left\{f_i: P \to P : i \in I_H\right\},\]
	such that 
	\begin{itemize}
		\item for all $i \in I_H$, $\hfl(f_i) < H$;
		\item if $i, j \in I_H$ satisfies that $i \neq j$, then $\hfl(f_i) \neq \hfl(f_j)$;
		\item $|I_H|$ is bounded by a function of $H$, increasing exponentially as $H$ increases. 
	\end{itemize}
	
	To construct a collection mentioned above, we use the notion of {\em Penner type symplectic automorphisms} that are mentioned above.
	We note that the $A_2$-Milnor fiber $P$ admits Penner type symplectic automorphisms, see Definition \ref{definition Penner type}. 
	We also recall that if $f$ is of Penner type, then $\hfl(f) = \hcat(f)$, as mentioned above.
	See Section \ref{subsection symplectic automorphisms of Penner type} for more details. 
	Moreover, on the $A_2$-Milnor fiber $P$, $\hcat(f)$ is computed in \cite{Bae-Choa-Jeong-Karabas-Lee25}.
	
	Combining all, we construct a collection $\{f_i: i \in I_H\}$ for all $H >0$, in Section \ref{section counting conjugacy classes in symplectic mapping class groups}. 
	The constructed collection satisfies that 
	\[|I_H| = e^H + O(1).\] 
	
	For the second part, we use the fact that if $\dim W = 4n \geq 6$, then every Dehn twist $\tau$ satisfies that $\tau^8$ is smoothly isotopic to the identity, as explained in \cite{Seidel14} and references therein. 
	Then, the proof of the first part, together with the above fact, proves the second half. 
\end{proof}

We also note that if there exists a conjugacy class $[f] \in C_{tri}(W)$ with $\hfl(f) >0$, $f$ posses some dynamical complexity that any symplectic homotopy cannot remove since $\hfl(f)$ is positive, but we can remove the dynamical complexity by smooth homotopy since $f$ is smoothly isotopic to the identity. 
Thus, it would be interesting if we can find a spectrum of different entropy values for $[f] \in C_{tri}(W)$. 
In the proof of Theorem \ref{thm application intro}, we can find a subset of the spectrum, and we can also show that the difference gap between two different entropy values tends to zero, as entropy values increase. 
See Corollary \ref{cor gap}.
%다음의 corollary를 Section 5에 집어 넣을 것.
%Cor: Under the same condition as Theorem \ref{thm application intro}, for every $\epsilon >0$, there exists $R_\epsilon >0$ such that if $r > R_\epsilon$, then there exists $f:W \to W$ satisfying that $|\hfl(f) -r| < \epsilon$.

\noindent{\bf Acknowledgment.}
During this work, the author has been supported by a KIAS Individual Grant (HP094401) at Korea Institute for Advanced Study.

\section{Setting, notation, and preliminaries}
\label{section setting, notation, and preliminaries}
In the section, we collect definitions, notation, preliminary results, etc. 
Section \ref{subsection symplectic topological background} specifies the geometric setting we will work on. 
Sections \ref{subsection Floer-theoretic entropy} and \ref{subsection categorial entropy} introduce Floer-theoretic and categorical entropies. 

\subsection{Symplectic topological background}
\label{subsection symplectic topological background}
We briefly summarize our geometric setting here. 
Especially, we introduce conditions of symplectic manifold $W$ and an automorphisms $f: W \to W$ should satisfy.
These conditions allow us to do Lagrangian Floer theory on $W$ and to study dynamics of $f$ using it.  
For more details, we refer the reader to \cite{Seidel08}.
Also, we define notation/terminologies that we will use in the later parts of the article.

In the current article, we work on Weinstein manifolds of finite type. 
By {\em finite type} Weinstein manifold $(W, d \lambda)$, we mean that $W$ can be constructed by attaching finitely many Weinstein handles up to Weinstein isotopy, or equivalently, $W$ admits a compact Lagrangian skeleton. 
We also assume that $W$ is {\em not} subcritical:
if $W$ is subcritical, its wrapped Fukaya category is the zero category, and the notion of categorical entropy becomes less interesting. 

Moreover, we would like to do Lagrangian Floer theory on Weinstein manifolds. 
Thus, we assume and define the following:
\begin{enumerate}[label=(\roman*), series=geometric_condition]
	\item First, we assume that $2c_1(W) =0$, thus various Fukaya categories of $W$ (such as wrapped, partially wrapped, etc.,) are $\mathbb{Z}$-graded. 
	\item There exists a fixed quadratic complex volume form $\eta^2$ on $W$. 
	It provides a squared phase map $\alpha: \mathrm{Gr}(TW) \to S^1$, where $\mathrm{Gr}(TW)$ means the Lagrangian Grassmannian. 
	\item Let $L$ be a Lagrangian submanifold of $W$. 
	Then, $L$ admits a natural section $s_L$ over $\mathrm{Gr}(W)$ sending $p \in L$ to $T_pL$. 
	We say that $L$ is {\em gradable} if $\alpha \circ s_L: L \to S^1$ has a lifting $\widetilde{\alpha \circ s_L}: L \to \mathbb{R}$. 
	Also, a pair of $L$ and a choice of lifting is called {\em graded Lagrangian}.
	\item We fix a {\em background class} $b \in H^2(W;\mathbb{Z}/2)$. 
	Then, we say a Lagrangian $L$ is {\em relatively pin} with respect to $b$, if the second Stiefel--Whitney class $w_2(L)$ and $b|_L$ coincide.  
\end{enumerate} 

By assuming (i)--(iv), we could define $\mathbb{Z}$-graded {\em wrapped} and {\em partially wrapped} Fukaya category of $W$. 
Let us define notation for them.
\begin{definition}
	\label{definition wrapped and partially wrapped Fukaya category}
	Let $W$ be a Weinstein manifold of finite type, satisfying (i)--(iv).
	\begin{enumerate}
		\item Let $\mathcal{W}(W)$ denote the wrapped Fukaya category of $W$. 
		\item Provided a stop $\Lambda \subset \partial_\infty W$, let $\mathcal{W}(W, \Lambda)$ denote the corresponding partially wrapped Fukaya category of $W$.
	\end{enumerate} 
\end{definition}
For more details on $\mathcal{W}(W)$ and $\mathcal{W}(W,\Lambda)$, we refer the reader to \cite{Sylvan19, Ganatra-Pardon-Shende20}

\begin{remark}
	\label{rmk grading issue}
	We note that if we work with $\mathbb{Z}/2$-graded Lagrangian Floer homologies, then we can simplify the above conditions (i)--(iv). 
\end{remark}

In the current article, we study the product of $W$ with itself, i.e., $W \times W$. 
Let us specify its Weinstein structure and define related notion.
\begin{definition}
	\label{definition product space}
	Let $W$ be a Weinstein manifold equipped with a Weinstein $1$-form $\lambda$. 
	\begin{enumerate}
		\item For simplicity, let $W^-$ denote the Weinstein manifold $(W, d (-\lambda))$, i.e., as a smooth manifold, it is the same as $W$, but its Liouville one-form has different sign from that of $W$. 
		\item By $W^- \times W$, we mean a Weinstein manifold 
		\[\left(W \times W, (p_1^*\left(-\lambda\right)) + p_2^*\lambda\right)\]
		where 
		\[p_1: W^- \times W \to W^-, (x,y) \mapsto x, \text{  and  } p_2:W^- \times W \to W, (x,y) \mapsto y.\]
		\item Let $\Delta$ denote the diagonal Lagrangian in $W^- \times W$, i.e., 
		\[\Delta:= \left\{(x,x) \in W^- \times W : x \in W\right\}.\]
	\end{enumerate}
\end{definition}

The following is trivial.
\begin{lem}
	\label{lem product space satisfies conditions}
	If $W$ is a Weinstein manifold satisfying (i)--(iv), so does $W^- \times W$. 
\end{lem}

From now on, we assume that every Weinstein manifold $W$ is of finite type and satisfies (i)--(iv) in the current article. 
We study compactly supported, symplectic automorphisms $f: W \to W$, which satisfy the following:
\begin{enumerate}[resume*=geometric_condition]
	\item $f^*\lambda - \lambda = dh$ for a compactly supported $h: W \to \mathbb{R}$, i.e., $f$ is an exact symplectic automorphism; 
	\item $f$ is equipped with a graded lift to the chosen Maslov cover;
	\item $f^*b =b$.
\end{enumerate}

If a compactly supported, symplectic automorphism $f:W \to W$ satisfies (v)--(vii), then $f$ induces an auto-equivalence on $\mathcal{W}(W)$ and also on $\mathcal{W}(W, \Lambda)$ for any $\Lambda$.
Objectwise, the induced functor sends an admissible Lagrangian $L$ to another admissible Lagrangian $f(L)$.
For more details, we refer the reader to \cite{Seidel08}. 
We set a notation for later.
\begin{definition}
	\label{definition notation phi}
	Let $f: W \to W$ be a compactly supported, symplectic automorphism satisfying (v)--(vii).
	We let 
	\[\phi: \mathcal{W}(W) \to \mathcal{W}(W), \quad \phi_\Lambda: \mathcal{W}(W;\Lambda) \to \mathcal{W}(W;\Lambda)\]
	denote the functors that $f$ induces on the corresponding categories.
\end{definition}
 
After having $\phi$, one can study the {\em categorical entropy} of it, (which will be defined in the next subsection,) and we can see it as an invariant of $f$.

If a compactly supported symplectic automorphism $f:W \to W$ is given, we can define another symplectic automorphism as follows:
\begin{definition}
	\label{definition product automorphism}
	Let $f: W \to W$ be a symplectic automorphism. 
	We let $F$ denote a symplectic automorphism defined on $W^- \times W$ as follows:
	\[F:W \times W \to W \times W, (x,y) \mapsto \left(x,f(y)\right).\]
\end{definition} 

We note that $F$ preserves the symplectic form, $p_1^*(-d\lambda) + p_2^*(d\lambda)$, i.e., $F$ is a symplectic automorphism.
Moreover, Lemma \ref{lem product automorphism} is trivial.
\begin{lem}
	\label{lem product automorphism}
	If $f:W \to W$ satisfies (vi) and (vii), then so does $F: W^- \times W \to W^- \times W$.
\end{lem}

Now, as similar to $f$, we would like to study the categorical counter part of $F$, i.e., the functors that $F$ induces on various Fukaya category. 
The problem is that even if $f$ satisfies (v), $F$ does not need to satisfy (v).
More precisely, if $f^*\lambda - \lambda = d h$ for a compactly supported $h: W \to \mathbb{R}$, then 
\[F^*\left[p_1^*\left((-\lambda )\right) + p_2^*\lambda\right] - \left[p_1^*\left((-\lambda )\right) + p_2^*\lambda\right] = d \left(h \circ p_2\right).\]
However, the support of $h \circ p_2$ is $W \times \left(\text{the support of  } h\right)$, which is not compact.
Thus, $F$ could not induce an {\em auto-equivalence} on $\mathcal{W}(W^- \times W)$ or its partially wrapped Fukaya categories.

Instead, $F$ induces a functor 
\[\widehat{F}: \mathcal{W}\left(W \times W, F^*\left[p_1^*(-\lambda)+p_2^*\lambda\right]\right) \to \mathcal{W}\left(W \times W, p_1^*(-\lambda) + p_2^*\lambda\right),\]
where $\mathcal{W}\left(W \times W; p_1^*(-\lambda) + p_2^*\lambda\right)$ and $\mathcal{W}\left(W \times W; F^*\left[p_1^*(-\lambda)+p^*_2\lambda\right]\right)$ denote the wrapped Fukaya categories of $W \times W$ when it admits a Liouville structure $p_1^*(-\lambda) + p_2^*\lambda$ and $F^*\left[p_1^*(-\lambda)+p_2\lambda\right]$, respectively.
We note that the difference between two wrapped Fukaya categories are controllable, since two Weinstein structures are Weinstein isotopic to each other by the following family of Liouville one-forms 
\begin{gather}
	\label{eqn isotopy family}
	\left\{\mu_s:=p_1^*\left((-\lambda )\right) + p_2^*\lambda + s d (h \circ p_2) : s \in [0,1]\right\}.
\end{gather}
We note that $h$ is compactly supported, thus the above one-parameter family of Weinstein structures gives us a Weinstein homotopy. 

Specifically, we define the following:
\begin{definition}
	\label{definition notation Phi}
	Let $(W, d\lambda)$ be a Weinstein manifold, $f: W \to W$ be a symplectic automorphism satisfying (v)--(vii), and $F = \mathbbm{1}_W \times f$ as defined in Definition \ref{definition product automorphism}.
	\begin{enumerate}
		\item Let $\varepsilon_0$ be the Weinstein homotopy 
		\item Let $\varepsilon$ denote the equivalence between two categories 
		\[\varepsilon: \mathcal{W}\left(W \times W, \mu_0\right) \to \mathcal{W}\left(W \times W, \mu_1\right).\]
		We note that the family $\{\mu_s\}$ is defined in Equation \eqref{eqn isotopy family} induces a deformation of Weinstein structure, and this deformation will induce $\varepsilon$ above. 
		\item Let $\Phi$ denote the following functor: 
		\[\Phi:= \widehat{F} \circ \varepsilon: \mathcal{W}\left(W \times W, \mu_0\right) \stackrel{\varepsilon}{\to} \mathcal{W}\left(W \times W, \mu_1\right) \stackrel{\widehat{F}}{\to} \mathcal{W}\left(W \times W, \mu_0\right). \]
		Or equivalently, $\Phi$ is the functor induced by a map $F \circ \varepsilon_0$. 
		\item Let $\Lambda$ be a stop on $\partial_\infty \left(W \times W\right)$, such that the deformation from $\{\mu_s\}$ does not change and such that $\Lambda$ is outside of the support of $F$. 
		Then, the deformation from $\{\mu_s\}$ and $F$ induce an auto-equivalence $\mathbf{\Phi_\Lambda}$ on the partially wrapped $\mathcal{W}\left(W \times W, \mu_0; \Lambda\right)$, i.e.,
		\[\Phi_\Lambda: \mathcal{W}\left(W \times W, \mu_0; \Lambda\right) \to \mathcal{W}\left(W \times W, \mu_0; \Lambda\right).\] 
	\end{enumerate} 
\end{definition}

The following definition would be not conventionally defined, but we use them for convenience. 
\begin{definition}
	\label{definition subdomain}
	For two Liouville manifolds $V$ and $W$, we say that $V$ is a {\em Liouville subdomain of $W$ up to completion}, if the following hold:
	$V$ and $W$ are the completions of two Liouville domains $V_0$ and $W_0$ such that $V_0$ is a Liouville subdomain of $W_0$.
	In other words, there exists an embedding 
	\[\iota: V_0 \to W_0 \text{  such that  } \iota^* \lambda_W = \lambda_V,\]
	where $\lambda_V$ and $\lambda_W$ are the Liouville one-forms on $V_0$ and $W_0$, respectively. 
\end{definition}
We note that if $V$ is a Liouville subdomain of $W$ up to completion, then there exists a cobordism $C$ such that $\partial_- C = \partial_\infty V$ and $\partial_+ C = \partial_\infty W$ and such that $W = V \cup C$. 
It would be a conventional way to describe the relation between $V$ and $W$, but for the sake of convenience, we use the phrase ``Liouville subdomain up to completion'', or simply ``Liouville subdomain'' if there is no chance of confusion.

\subsection{Floer-theoretic entropy}
\label{subsection Floer-theoretic entropy}
The purpose of Section \ref{subsection Floer-theoretic entropy} is to define the notion of {\em Floer-theoretic} entropy which is given in Equation \eqref{eqn Floer theoretic entropy}, as well as to explain its basic properties that we will use later. 

We will work in the setting explained in Section \ref{subsection symplectic topological background}.
To the define Floer-theoretic entropy, we need to introduce the {\em fixed-point Floer homology} of $f:W \to W$. 
However, we omit the introduction in the current paper, since there are many of references defining the notion, including \cite{Uljarevic17, McLean12}. 

\begin{definition}
	\label{definition Floer-theoretic entropy}
	Let $W$ be a Liouville manifold, and let $f:W \to W$ be a compactly supported, exact symplectic automorphism. 
	The {\em Floer-theoretic entropy} of $f$ is given as 
	\[\hfl(f) := \limsup_{n \to \infty} \frac{1}{n} \log_+ \dim \HF_{fixed}(f^n),\]
	where $\HF_{fixed}(f^n)$ means the fixed-point Floer homology of $f^n$ and $\log_+$ is the function defined as follows:
	\[\log_+: \mathbb{Z}_{\geq 0} \to \mathbb{R} \text{  such that  } \log_+(0) = 0, \log_+(k) = \log k \text{  for  } k \in \mathbb{N}.\]
\end{definition}

\begin{remark}
	\label{rmk Floer-theoretic entropy}
	\mbox{}
	\begin{enumerate}
		\item Compared to Equation \eqref{eqn Floer theoretic entropy}, we are using $\log_+$ in Definition \ref{definition Floer-theoretic entropy}.
		The purpose of using $\log_+$ is that it gives us a well-defined, non-negative Floer-theoretic entropy, even if $\HF_{fixed}(f^n) = 0$ for infinitely many $n$. 
		\item We note that if $f, g \in \pi_0 \left(\textit{Symp}(W)\right)$ are symplectically isotopic to each other, then $\hfl(f) = \hfl(g)$, since they have the same fixed-point Floer homology by \cite[Theorem 2.43]{Uljarevic17}.
		Moreover, if $f, g \in \pi_0 \left(\textit{Symp}(W)\right)$ are conjugate to each other, then they also have the same Floer-theoretic entropy, as stated in \cite{Smith12}.
	\end{enumerate}
\end{remark}

We state a {\em local-to-global} property of Floer-theoretic entropy, which we will use in Section \ref{section counting conjugacy classes in symplectic mapping class groups}.
The property is proven in \cite{Kim-Kwon26}.
\begin{thm}[Theorem 3.2 of \cite{Kim-Kwon26}]
	\label{thm local to global property}
	Let $V_0, W_0$ be Liouville domains such that $V_0$ is a Liouville subdomain of $W_0$. 
	Let $f:V_0 \to V_0$ be an exact symplectic automorphism whose support is contained in the interior of $V_0$. 
	Then, $f$ can extend to an exact symplectic automorphism $\widetilde{f}:W_0 \to W_0$ by setting 
	\[\widetilde{f}(x) = \begin{cases}
		f(x), &\text{  if  } x \in V_0, \\
		x, & \text{  if  } x \in W_0 \setminus V_0.
	\end{cases}\]
	Then, $\hfl(f) = \hfl(\widetilde{f})$.
\end{thm}

\subsection{Categorical entropy}
\label{subsection categorial entropy}
In Section \ref{subsection symplectic topological background}, we introduced our setting. 
Especially, our given Weinstein manifold $W$ and a symplectic automorphism $f: W \to W$ are assume to satisfy conditions (i)--(vii).
Then, the given dynamical system induces a category, a (partially) wrapped Fukaya category, and a functor, denoted $\phi$ or $\phi_\Lambda$ in Definition \ref{definition notation phi}, on it.
Now, we can say that the categorical invariants of $\phi$ (or $\phi_\Lambda$) are symplectic topological invariants of $f$, defined up to Weinstein homotopy. 
In the current subsection, we define and explain known results about a categorical invariant that is called {\em categorical entropy}.

Originally, the notion of categorical entropy is defined in \cite{Dimitrov-Haiden-Katzarkov-Kontsevich}, for a triangulated category with a split-generator. 
In the current article, we will use another definition introduced in \cite{Bae-Choa-Jeong-Karabas-Lee25}, which is equivalent to the original definition by \cite[Proposition 2.13]{Bae-Choa-Jeong-Karabas-Lee25}.
The other definition is defined for a pretriangulated category admitting a finite collection of generators.

\begin{definition}
	\label{definition categorical entropy}
	Let $\mathcal{C}$ be a pretriangulated, nonzero $A_\infty$ (or dg) category.
	Let us assume that $\mathcal{C}$ admits a finite collection of generators $G_1, \dots, G_m$, i.e., for any object $E \in \mathcal{C}$, $E$ can be written as a twisted complex consisting of (multiple copies of) shifted $G_1, \dots, G_m$.
	In other words, $E$ should be equivalent to a twisted complex $\left(E_i, f_{i,j}\right)$ such that $E_i$ is a one of $\{G_1[d], \dots, G_m[d] : d \in \mathbb{Z}\}$. 
	\begin{enumerate}
		\item For an object $E \in \mathcal{C}$ the {\em complexity of $E$ with respect to $G_1, \dots, G_m$} is defined as 
		\[\delta(G_1,\dots,G_m;E) := \inf \left\{\text{  the number of  }G_i[d] \text{  in a twisted complex  } K : K \text{  is equivalent to  } E \oplus E' \text{  for some  } E' \in \mathcal{C}\right\}.\]
		\item Let $\psi:\mathcal{C} \to \mathcal{C}$ be an auto-equivalence. 
		Then, the {\em categorical entropy of $\psi$}, denoted as $\mathbf{\hcat(\psi)}$ is defined as follows:
		\[h(\psi) := \lim_{n\to \infty} \tfrac{1}{n} \log \delta\left(G_1, \dots, G_m;\psi^n(G_1 \oplus \dots \oplus G_m)\right).\] 
	\end{enumerate} 
\end{definition}

\begin{remark}
	\label{rmk at t}
	We note that the original definition of categorical entropy requires to fix a number $t \in \mathbb{R}$, and define it at $t$. 
	Roughly speaking, the number $t$ is used to measure how much shift is involved in the dynamics of an auto-equivalence. 
	However, in the current article, we only need to use the categorical entropy at $t=0$. 
	Thus, in Definition \ref{definition categorical entropy}, we omitted the factor $t$ for simplicity.
	Moreover, the original definition of categorical entropy is defined for an arbitrary functor that is not necessarily to be an auto-equivalence, and also the category is not necessarily to be nonzero. 
	However, we restricted to our attention to auto-equivalences on a nonzero category.
	Since the main geometric setting of the paper is non-subcritical Weinstein manifolds and symplectic automorphisms on them, our restricted interests can be justified.  
	By doing so, we can assume that the categorical entropy cannot be $-\infty$, which is allowed in \cite{Dimitrov-Haiden-Katzarkov-Kontsevich}.
	Our assumption can be justified because our setting is based on symplectic topology.
	Thus, the categories and functors we are interested in would satisfy the assumptions. 
\end{remark}

By \cite{Chantraine-Dimitroglou-Ghiggini-Golovko24, Ganatra-Pardon-Shende24}, it is known that if $W$ is a Weinstein manifold of finite-type, the wrapped Fukaya category of it admits a finite collection of Lagrangian disks that generate the wrapped Fukaya category. 
Thus, one can define the categorical entropy of a symplectic automorphism $f: W \to W$ as the categorical entropy of the induced functor.

One problem is that there exist various versions of Fukaya categories. 
For example, for any stop $\Lambda$ of a Weinstein manifold $W$, there exists the notion of partially wrapped Fukaya category $\mathcal{W}(W;\Lambda)$.
And, \cite{Chantraine-Dimitroglou-Ghiggini-Golovko24, Ganatra-Pardon-Shende24} showed that $\mathcal{W}(W;\Lambda)$ has generators.
Thus, the notion of categorical entropy is defined well on each partially wrapped Fukaya category.
 
A natural question asks about the relation between categorical entropies on wrapped and partially wrapped Fukaya categories.
There is a known answer for the question: 
\cite{Bae-Choa-Jeong-Karabas-Lee25} showed that if $f:W \to W$ is compactly supported, then the induced functors on wrapped and partially wrapped Fukaya categories have the same categorical entropy.
Remark \ref{rmk key idea} explains the key idea that we will use later. 
\begin{remark}
	\label{rmk key idea}
	The key idea of the result given in \cite{Bae-Choa-Jeong-Karabas-Lee25} is that if $f$ is compactly supported, then near any chosen stop $\Lambda$, $f$ acts as the identity automorphism.
	Thus, the induced functor has no nontrivial effect on {\em linking disks} that are generators associated to the choice of $\Lambda$.
	It implies that the dynamical complexity of the induced functor can be computed without considering the generators associated to $\Lambda$.
	It leads us to the equality between categorical entropies on various (partially) wrapped Fukaya categories. 
\end{remark} 

Than, we can define the categorical entropy of a compactly supported symplectic automorphisms $f:W \to W$ as follows:
\begin{definition}
	\label{definition categorical entropy of symplectic dynamics}
	Let $W$ be a Weinstein manifold of finite type and let $f:W \to W$ be a compactly supported symplectic automorphism.
	Then, the {\em categorical entropy of $f$} is defined to be that of the functor $\phi: \mathcal{W}(W) \to \mathcal{W}(W)$ that $f$ induces. 
\end{definition}

\section{Proof of Theorem \ref{thm intro}}
\label{section proof of the main theorem}
Let $W$ be a Weinstein manifold of finite type, having nonzero wrapped Fukaya category $\mathcal{W}(W)$. 
We fix a compactly supported symplectic automorphism $f: W \to W$ that induces functors on wrapped and various partially wrapped Fukaya category. 
In the section, we prove our main result, Theorem \ref{thm intro}, by employing notation defined in Section \ref{section setting, notation, and preliminaries}.

Our proof consists of the following three steps.
\begin{itemize}
	\item(Section \ref{subsection the first step}) The first step is to translate $\hfl(f)$ in terms of a partially wrapped Fukaya category. 
	More specifically, we will fix a specific stop $\Lambda_0$ and we will show that 
	\begin{gather}
		\label{eqn the first step}
		\hfl(f) \leq \hcat(\Phi_{\Lambda_0}).
	\end{gather}  
	\item(Section \ref{subsection the second step}) The second step is to show that the fixed stop $\Lambda_0$ does not affect on the computation of categorical entropies of $F$. 
	More specifically, we show 
	\begin{gather}
		\label{eqn the second step}
		\hcat(\Phi_{\Lambda_0}) = \hcat(\Phi).
	\end{gather}
	\item(Section \ref{subsection the third step}) The final step is to show that two functors that are induced by $f$ and $F = \mathbbm{1}_W \times f$ respectively satisfy the following inequality:
	\begin{gather}
		\label{eqn the third step}
		\hcat(\Phi) \leq \hcat(\phi) =\hcat(f).
	\end{gather}
\end{itemize} 

Combining Equations \eqref{eqn the first step}--\eqref{eqn the third step} proves 
\[\hfl(f) \leq \hcat(f),\]
i.e., Theorem \ref{thm intro}.

\subsection{The first step}
\label{subsection the first step}
In this subsection, we show that the Floer-theoretic entropy of $f$ is bounded above by the categorical entropy of $\Phi_\Lambda$, defined on a partially wrapped Fukaya category of $W^- \times W$. 

First, we simply use \cite[Theorem 2.6]{Kim-Kwon26}, stated below.
To state \cite[Theorem 2.6]{Kim-Kwon26}, we need to define the following notation: 
\cite[Theorem 2.6]{Kim-Kwon26} relates $\HF_{\mathrm{fixed}}(f)$ with a Lagrangian Floer homology in $W^- \times W$, between $\Delta$ and the graph of $f$, i.e., 
\[\left\{(x,f(x)) : x \in W\right\} = F(\Delta).\]
A problem is that $\Delta$ and $F(\Delta)$ do not intersect transversally:
Actually, in a cylindrical part of $W^- \times W$, $\Delta$ and $F(\Delta)$ coincide since $f$ is compactly supported. 
To handle this problem, \cite{Kim-Kwon26} took a small deformation of $\Delta$, which is given as an arbitrary small positive Reeb flow in a cylindrical end of $W^- \times W$ as well as a compactly supported deformation.
More precisely, \cite{Kim-Kwon26} fixes a Hamiltonian $H$ that is linear near the boundary with positive slope $\sigma$ such that $\sigma$ is strictly less that the minimum of the lengths of Reeb chords from $\partial_\infty \Delta$ to $\partial_\infty F(\Delta) = \partial_\infty \Delta$.
Then, one can take the perturbation of $\Delta$ by time-$1$ flow of the Hamiltonian vector field of $H$, and it gives us the usual Lagrangian Floer homology of two transversally intersecting Lagrangians, i.e., the perturbed $\Delta$ and $F(\Delta)$. 
Let $\HF_{\mathrm{small}}\left(\Delta, F(\Delta)\right)$ denote the corresponding Lagrangian Floer homology.  

\begin{thm}[Theorem 2.6 of \cite{Kim-Kwon26}]
	\label{thm folklore}
	For every compactly supported $f: W \to W$, there exists a canonical isomorphism
	\[\HF_{\mathrm{fixed}}(f) \simeq \HF_{\mathrm{small}}\left(\Delta, F(\Delta)\right),\]
	where $\HF_{\mathrm{fixed}}(f)$ is the fixed-point Floer homology of $f$.
\end{thm}

Now, let us assume that the $\sigma$ (the slope of $H$ that we used to define $\HF_{\mathrm{samll}}(\Delta, F(\Delta))$ above) is strictly smaller that the {\em half} of the minimum of the lengths of Reeb chords from $\partial_\infty \Delta$ to $\partial_\infty F(\Delta)$.
And, we set $\Lambda_0$ be the Legendrian in $\partial_\infty \left(W^- \times W\right)$ obtained by taking time-$\left(\tfrac{3}{2}\sigma\right)$ Reeb flow of another Legendrian $\partial_\infty \Delta$. 
Then, we consider the partially wrapped Fukaya category of $\left(W^- \times W, \Lambda_0\right)$. 

By definition of the partially wrapped Fukaya category, $\HF_{\mathrm{small}}(\Delta, F(\Delta))$ coincides with the morphism space from $\Delta$ to $F(\Delta)$ in $\mathcal{W}\left(W^- \times W; \Lambda_0\right)$.
In other words, we have
\begin{gather}
	\label{eqn first step 1}
	\HF_{\mathrm{fixed}}(f) \simeq \HF_{\mathrm{small}}(\Delta, F(\Delta)) = \HF_{\left(W^- \times W; \Lambda_0\right)}\left(\Delta, F(\Delta)\right),
\end{gather}
where the last one is the morphism space in (the derived) wrapped Fukaya category $\mathcal{W}(W^- \times W; \Lambda_0)$. 

\begin{figure}[h]
	\centering
	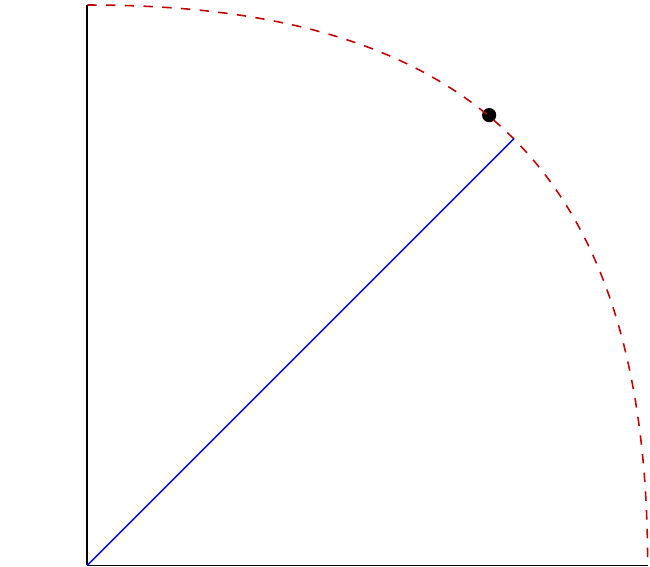		
	\caption{The horizontal black line on the bottom describes $W^-$ and the vertical black line on the left describes $W$. 
		The blue line corresponds to the diagonal Lagrangian $\Delta$, and the red dashed curve represents $\partial_\infty \left(W^- \times W\right)$. 
		The black dot is the stop $\Lambda_0$.
		On the left vertical line, we have a compact subset of $W$, which is the support of $f$.
		Thus, the support of $F$ is the region below of the black dashed line. 
		Since $\Delta_0$ is outside of the support of $F$, it does not change by the deformation $\varepsilon_0$.}
	\label{figure first step}
\end{figure}

Moreover, we can easily check that $\Lambda_0$ does not change by the deformation $\varepsilon_0$ defined in Definition \ref{definition notation Phi}.
Figure \ref{figure first step} explains it pictorially. 
Thus, we have an auto-equivalence 
\[\Phi_{\Lambda_0}: \mathcal{W}(W^- \times W; \Lambda_0) \to \mathcal{W}(W^- \times W; \Lambda_0).\]  

Since $F(\Delta)$ is a Lagrangian representing an object $\Phi_{\Lambda_0}(\Delta)$, Equation \eqref{eqn first step 1} becomes 
\[\HF_{\mathrm{fixed}}(f) \simeq \HF_{\left(W^- \times W; \Lambda_0\right)}\left(\Delta, \Phi_{\Lambda_0}(\Delta)\right).\]
Also, by the same logic, we can put $f^n$ in the position of $f$, and it leads us 
\begin{align}
	\label{eqn first step 2}
	\begin{split}
		\hfl(f) &= \limsup_{n \to \infty} \tfrac{1}{n} \log_+ \dim \HF_{\mathrm{fixed}}(f^n) \\
	&= \limsup_{n \to \infty} \tfrac{1}{n} \log_+ \dim \HF_{\left(W^- \times W; \Lambda_0\right)}\left(\Delta, \Phi_{\Lambda_0}^n(\Delta)\right).
	\end{split}
\end{align}

Now, we observe the following facts:
\begin{enumerate}
	\item By \cite{Chantraine-Dimitroglou-Ghiggini-Golovko24, Ganatra-Pardon-Shende24}, there exists a generator $G$ that is given as a disjoint, finite union of Lagrangian cocores and linking disks associated to $\Lambda_0$.
	A detailed description of $G$ will be given in the next subsection, and we use only the existence of a generator $G$, represented by Lagrangians, in the current subsection.
	\item By definition of partially wrapped Fukaya category, $\HF_{\left(W^- \times W; \Lambda_0\right)}\left(\Delta, G\right)$ has to be finite-dimensional. 
	\item Since $G$ is a generator, we can find an iterated cone of multiple copies of $G$ which is equivalent to $\Phi^n_{\Lambda_0}(\Delta)$. 
	Let $K$ denote such an iterated cone, and let $\ell$ denote the number of shifted $G$ in $K$. 
	Then, one can easily see that 
	\[\dim \HF_{\left(W^- \times W; \Lambda_0\right)}(\Delta, \Phi^n_{\Lambda_0}\Delta) \leq \ell \cdot \dim \HF_{\left(W^- \times W; \Lambda_0\right)}(\Delta, G).\]
	Using the complexity function $\delta$ defined in Definition \ref{definition categorical entropy}, the above inequality implies that 
	\[\dim \HF_{\left(W^- \times W; \Lambda_0\right)}(\Delta, \Phi^n_{\Lambda_0}\Delta) \leq \delta(G, \Phi^n_{\Lambda_0}\Delta) \cdot \dim \HF_{\left(W^- \times W; \Lambda_0\right)}(\Delta, G).\]
	\item From definition of the complexity $\delta$, it is trivial to see that
	\[\delta(G,\Phi^n_{\Lambda_0}\Delta) \leq \delta(G,\Phi^n_{\Lambda_0}G) \cdot \delta(\Phi^n_{\Lambda_0}G, \Phi^n_{\Lambda_0}\Delta) = \delta(G,\Phi^n_{\Lambda_0}G) \cdot \delta(G, \Delta). \]
	For the last equality, we used that $\Phi_{\Lambda_0}$ is an auto-equivalence. 
\end{enumerate}

Finally, combining (3) and (4), we have
\[\dim \HF_{\left(W^- \times W; \Lambda_0\right)}(\Delta, \Phi^n_{\Lambda_0}\Delta) \leq \delta(G,\Phi^n_{\Lambda_0}G) \cdot \delta(G, \Delta) \cdot \dim \HF_{\left(W^- \times W; \Lambda_0\right)}(\Delta, G).\]
Since the exponential growth of the left-handed side is $\hfl(f)$ and that of the right-handed side is $\hcat(\Phi_{\Lambda_0})$,  
we have Equation \eqref{eqn the first step}, i.e., $\hfl(f) \leq \hcat(\Phi_{\Lambda_0})$.

\subsection{The second step}
\label{subsection the second step}
The goal of the subsection is to show that Equation \eqref{eqn the second step} holds. 
In other words, we would like to show that two auto-equivalences $\Phi_{\Lambda_0}$ and $\Phi$, defined on different wrapped Fukaya categories, have the same categorical entropy.
To show that, we specify the generator $G$ of the partial wrapped Fukaya category $\mathcal{W}(W^- \times W; \Lambda_0)$, which is used at the last part of Section \ref{subsection the first step}.

We recall that $W$ is a Weinstein manifold of finite type. 
It implies that by \cite{Chantraine-Dimitroglou-Ghiggini-Golovko24, Ganatra-Pardon-Shende24}, there exits a finite collection of pairwise disjoint Lagrangian disks $G_1, \dots, G_m$ such that $\{G_1, \dots, G_m\}$ generates the wrapped Fukaya category $\mathcal{W}(W)$.
Also, the same collection generates $\mathcal{W}(W^-)$, i.e., the wrapped Fukaya category of the same manifold having the opposite-signed Weinstein structure.

We apply the K\"{u}nneth formula from \cite[Section 8]{Ganatra-Pardon-Shende24}.
Then, a collection of Lagrangians 
\[\left\{G_i \times G_j \subset W^- \times W : i, j = 1, \dots, m\right\} \]
generates $\mathcal{W}(W^- \times W)$. 
Additionally, when we fix a specific stop $\Lambda_0 \subset \partial_\infty \left(W^- \times W\right)$, there exists a finite collection of linking disks $D_1, \dots, D_\ell$ such that 
\[\left\{D_k, G_i \times G_j \subset W^- \times W: k= 1, \dots, \ell, i, j = 1, \dots, m\right\}\]
generates $\mathcal{W}(W^- \times W; \Lambda_0)$. 
We note that linking disks $D_1, \dots, D_\ell$ are contained in a cylindrical end of $W^- \times W$, near $\Lambda_0$.
Thus, one can choose them away from $G_{i,j}$. 
We remark that the generator $G$, used in the last part of Section \ref{subsection the first step}, is the disjoint union of $G_i \times G_j$ and $D_k$, i.e., 
\[G = \left(\bigoplus_{i,j} \left(G_i \times G_j\right)\right) \oplus \left(\bigoplus_{k} D_k\right).\] 

From the definition of $D_k$, one can easily see that the symplectic automorphism $F$ and the deformation$\varepsilon_0$ do not change $D_k$. 
Thus,
\begin{gather}
	\label{eqn second step}
	\Phi_{\Lambda_0}(D_k) = D_k \text{  for all  } k =1, \dots, \ell.
\end{gather}
It implies that if $\mathcal{D}$ is a thick subcategory of $\mathcal{W}(W^- \times W; \Lambda_0)$ generated by $D_1, \dots, D_\ell$, then the restriction of $\Phi_{\Lambda_0}$ defines an auto-equivalence $\Phi_{\Lambda_0}|_{\mathcal{D}}: \mathcal{D} \to \mathcal{D}$. 
Moreover, the restriction has the zero categorical entropy, because of Equation \eqref{eqn second step}. 

We apply \cite[Theorem 1.1]{Bae-Choa-Jeong-Karabas-Lee25} stated as follows:
\begin{thm}[Theorem 1.1 of \cite{Bae-Choa-Jeong-Karabas-Lee25}]
	\label{thm stop removal}
	Let $\mathcal{C}$ be a pretriangulated $A_\infty$ category with generators, and let $\mathcal{D}$ be a thick subcategory of $\mathcal{C}$ with generators. 
	Let $\Psi_\mathcal{C}: \mathcal{C} \to \mathcal{C}$ be a functor preserving $\mathcal{D}$ so that there are induced functors 
	\[\Psi_\mathcal{D}: \mathcal{D} \to \mathcal{D} \text{  and  } \Psi_{\mathcal{C}/\mathcal{D}}: \mathcal{C}/\mathcal{D} \to \mathcal{C}/\mathcal{D}.\]
	Then, the following equality holds:
	\[\hcat(\Psi_{\mathcal{C}/\mathcal{D}}) \leq \hcat(\Psi_\mathcal{C}) \leq \max \left\{\hcat(\Psi_{\mathcal{C}/\mathcal{D}}), \hcat(\Psi_\mathcal{D})\right\}.\]
\end{thm}

In our case, $\mathcal{C}$ is the partially wrapped Fukaya category $\mathcal{W}(W^- \times W; \Lambda_0)$, $\mathcal{D}$ is the thick subcategory generated by $D_k$.
Then, the quotient 
\[\mathcal{C}/\mathcal{D} = \mathcal{W}(W^- \times W; \Lambda_0) / \mathcal{D} = \mathcal{W}(W^- \times W).\]
And for the functors, $\Psi_\mathcal{C} = \Phi_{\Lambda_0}$ and $\Psi_{\mathcal{C}/\mathcal{D}} = \Phi$. 

Applying Theorem \ref{thm stop removal} provides us 
\[\hcat(\Phi) \leq \hcat(\Phi_{\Lambda_0}) \leq \max \left\{ \hcat(\Phi), \hcat(\Phi_{\Lambda_0}|_\mathcal{D})=0\right\}.\]
As mentioned in Remark \ref{rmk at t}, our categorical entropy is always a non-negative number, so the above implies that 
\[\hcat(\Phi_{\Lambda_0}) = \hcat(\Phi),\]
which is the goal of this subsection, i.e., Equation \eqref{eqn the second step}.

\subsection{The third step}
\label{subsection the third step}
The final step, i.e., showing
\[\hcat(\Phi) \leq \hcat(f),\]
is the simplest step. 
To prove the simplest step, we need the following equivalence 
\[\mathcal{T}: \mathcal{W}(W^- \times W) \stackrel{\simeq}{\to} \mathcal{W}(W^-) \otimes \mathcal{W}(W),\]
provided in \cite[Section 8]{Ganatra-Pardon-Shende24}.

To be more precise, let us recall that $\hcat(\Phi)$ is defined to be the exponential growth ratio of 
\[\delta\left(\left\{G_i \times G_j\right\}, \Phi^n\left(\oplus_{i,j} \left(G_i \times G_j\right)\right)\right) \leq \sum_{(i_0,j_0)}\delta\left(\left\{G_i \times G_j\right\}, \Phi^n\left(G_{i_0} \times G_{j_0}\right)\right),\]
as $n \to \infty$.
To compute it, let us fix $i_0$ and $j_0$. 
Then, we have 
\[\Phi^n\left(G_{i_0} \times G_{j_0}\right) = G_{i_0} \times f^n\left(G_{j_0}\right).\]
We apply $\mathcal{T}$, 
\[\mathcal{T} \left(\Phi^n\left(G_{i_0} \times G_{j_0}\right)\right) \simeq G_{i_0} \otimes f^n\left(G_{j_0}\right).\]

We note that tensoring with $G_{i_0}$ is an exact functor, i.e., 
\[G_{i_0} \otimes \text{--} : \mathcal{W}(W) \to \mathcal{W}(W^-) \otimes \mathcal{W}(W)\]
is an exact functor.
For more details, see \cite{Ornaghi25}.
Thus, if we have an iterated cone representing $f^n\left(G_{j_0}\right)$, it induces an iterated cone representing $G_{i_0} \otimes f^n\left(G_{j_0}\right)$, and two iterated cones have the same length.
It implies that 
\[\delta\left(\{G_i \times G_j\}, \Phi^n\left(G_{i_0} \times G_{j_0}\right)\right) \leq \delta\left(\left\{G_i\right\}, f^n\left(G_{j_0}\right)\right).\]
By measuring the exponential growth ratio, we can prove the desired inequality, i.e., Equation \eqref{eqn the third step}.

\section{The sharpness of the inequality in Theorem \ref{thm intro}}
\label{section the sharpness of the inequality}
A natural mathematical question following Theorem \ref{thm intro} asks the {\em sharpness} of the inequality 
\[\hfl(f) \leq \hcat(f).\]
In this section, we provide examples of two types: the first type example gives a strict inequality (see Section \ref{subsection the torus example}) and the second gives an equality (see Section \ref{subsection symplectic automorphisms of Penner type}).

\subsection{The torus example}
\label{subsection the torus example}
In the current subsection, we give an example of a symplectic automorphism $f$ such that 
\[\hfl(f)=0 < \hcat(f).\]
However, this example automorphism is defined on a {\em closed} symplectic manifold, differently from all other parts of the article. 
The purpose of considering this exceptional setting is to explain the difference between Floer-theoretic and categorical entropies, at least on a philosophical level. 

The closed symplectic manifold which we consider is $M:=\mathbb{T} \times \mathbb{T}$, where $\mathbb{T}$ is the $2$-dimensional torus equipped with a symplectic form. 
Thus, $M$ is a product of two symplectic manifolds, and it defines a symplectic structure of $M$. 

In $\mathbb{T}$, we can fix a meridian and longitude, which we will denote by $m$ and $\ell$. 
Then, we know that 
\[\left\{m \times m, m \times \ell, \ell \times m, \ell \times \ell\right\}\]
generate the Fukaya category of $M$, as stated in \cite{Abouzaid-Smith10}.
Thus, the categorical entropy of an auto-equivalence on the Fukaya category is well-defined. 

We consider an automorphism on $\mathbb{T}$ given as a matrix 
\[A = \begin{bmatrix}
	2 & 1 \\ 1 & 1 
\end{bmatrix}. \]
Here, the first (resp.\ second) coordinate corresponds to the meridian (resp.\ longitude) direction. 
And we choose a irrational real number $\alpha \in \mathbb{R}\setminus \mathbb{Q}$, and fix a translation map $\tau_\alpha: \mathbb{T} \to \mathbb{T}$ which is defined as 
\[\tau_\alpha(x,y) = (x + \alpha, y).\]

Since $M$ is product of two copies of $\mathbb{T}$, the following gives a symplectic automorphism defined on $M$: 
\[f: M \to M, f(p,q) = \left(A p, \tau_\alpha(q)\right).\]
One can easily check that any power of $f$, $f^n$, do not have any fixed-point, because of the translation factor. 
Thus, the fixed-point Floer homology of $f^n$ should vanish, and $\hfl(f) = 0$. 

To show that $\hcat(f) >0$, we recall that in Section \ref{subsection the first step}, we proved that the exponential growth of 
\[\dim \HF^*\left(\Delta, F^n\left(\Delta\right)\right)\]
bounds from below the categorical entropy of $F$. 
The same proof can show that the exponential growth of 
\[\dim \HF^*\left(m \times m, f^n\left(m \times m\right)\right)\] 
is smaller than or equal to $\hcat(f)$. 
And, one can easily compute that the above exponential growth is the log of the greatest Eigen-value of $A$, i.e., $\tfrac{3 + \sqrt{5}}{2}$. 
Thus, $\hcat(f) > 0$.

\subsection{Symplectic automorphisms of Penner type}
\label{subsection symplectic automorphisms of Penner type}
We introduce a construction of symplectic automorphisms whose categorical and Floer-theoretic entropy coincide. 

These symplectic automorphisms are called {\em Penner type}.
They are originally defined in \cite{Lee24}, and \cite{Bae-Choa-Jeong-Karabas-Lee25} gives an algorithm computing their categorical entropy.
Recently, \cite{Kim-Kwon26} computed a lower-bound of their Floer-theoretic entropy.
Interestingly, the lower-bounds coincide with the categorical entropy of them, which can be computed by \cite{Bae-Choa-Jeong-Karabas-Lee25}.
Shortly, 
\[\text{a lower-bound by \cite{Kim-Kwon26}  } \leq \text{  Floer-theoretic entropy  } \leq \text{  categorical entropy  } = \text{  the lower-bound by \cite{Kim-Kwon26}}. \]

In order to define symplectic automorphisms of Penner type and to provide more details, we need to introduce some definition/notation first. 

\subsubsection{Plumbing spaces}
\label{subsubsection plumbing spaces}
We introduce the {\em plumbing spaces} on which symplectic automorphisms will be defined.
In the current paper, we briefly introduce the construction of plumbing spaces.
For an eager reader, we refer the reader to \cite{Abouzaid-Smith12, Karabas-Lee25}.

For convenience, let us fix a natural number $n \in \mathbb{N}$. 
We first introduce the simplest case, a plumbing of two $T^*S^n$. 
Simply speaking, a plumbing space is obtained by identifying a small region of $T^*S^n$ and another small region of the other $T^*S^n$, by flipping the base and fiber directions.  

More rigorously, we plumb two disk cotangent bundles $D^*S^n$ at one point. 
To do that, let $S_i = S^n$ for $i= 1, 2$.
Then, we choose plumbing points $p \in S_1$ and $q \in S_2$. 
Since $S_i$ is a smooth manifold, there exists a small neighborhood $U_p \subset S_1$ (resp.\ $V_q \subset S_2$) such that 
\[U_p \simeq \mathbb{D}^n \simeq V_q.\]
And, moreover, $D^*U_p \subset D^*S_1$ and $D^*V_q \subset D^*S_2$. 

We note that 
\[D^*U_p \simeq D^*\mathbb{D}^n = \mathbb{D}^n \times \mathbb{D}^n \simeq D^*V_q.\]
Using this, one can construct a map
\begin{gather}
	\label{eqn identification map}
	\psi: D^*U_p \simeq \mathbb{D}^n \times \mathbb{D}^n \to \mathbb{D}^n \times \mathbb{D}^n \simeq D^*V_q, \quad (y,x) \mapsto (x,-y).
\end{gather}
Since $\psi$ preserves the standard symplectic form $d y \wedge d x$, thus $\psi$ is a symplectomorphism. 

Now, we can plumb $D^*S_1$ and $D^*S_2$.
The resulting plumbing space is defined to be a convex completion of a smoothing of the following
\[D^*S_1 \sqcup D^*S_2 / \sim, \text{  where for all  } (y,x) \in D^*U_p, (y,x) \sim \psi(y,x). \]
Since $\psi$ is a symplectomorphism, the constructed space has a natural symplectic structure that is inherited from those of $T^*S^n$s. 
Moreover, the Weinstein structures on $T^*S^n$s induce a Weinstein structure on the plumbing space. 
For more details on the Weinstein structure, see \cite{Abouzaid-Smith12} for a classical explanation, and see \cite{Karabas-Lee25} for one using the notion of Weinstein sectorial covering. 

Note that the above plumbing procedure is a {\em local procedure}: the whole procedure identifies a small region of $T*S^n$, and outside of the small regions, we did not do anything. 
Thus, one can perform multiple plumbing procedures, connecting multiple copies of $T^*S^n$, simultaneously.
A {\em plumbing space} means a Weinstein manifold constructed by plumbing multiple cotangent bundles. 

In \cite{Karabas-Lee25}, the authors used a quiver $Q$ describe a plumbing space. 
\begin{definition}
	\label{definition plumbing space}
	Let $Q$ be a finite quiver, i.e., a directed graph having finitely many vertices and arrows. 
	Let $P_n(Q)$ denote the plumbing space constructed as follows: 
	\begin{itemize}
		\item We first assign one $T^*S^n$ to each vertex of $Q$, and let $T^*S_v$ denote the one assigned to a vertex $v$ of $Q$;
		\item If there is an arrow $e = v \to w$ connecting two vertices $v$ and $w$, the we plumb $T^*S_v$ and $T^*S_w$;
		\item The direction of $e = v \to w$ determines the direction of the identifying map $\phi$ (given in Equation \eqref{eqn identification map}), i.e., a small region of $T^*S_v$ (resp.\ $T^*S_w$) becomes the domain (resp.\ target) of $\phi$.
	\end{itemize} 
	We say that $Q$ is the {\em plumbing pattern} of $P_n(Q)$. 
\end{definition}

\begin{remark}
	\label{rmk a general plumbing}
	\mbox{}
	\begin{enumerate}
		\item The notion of plumbing space in \cite{Karabas-Lee25} is more general than Definition \ref{definition plumbing space} in the following two ways: 
		First, \cite{Karabas-Lee25} allows us to plumb cotangent bundles that are different from $T^*S^n$. 
		We note that the plumbing procedure does not use any topology of $T^*S^n$, thus one can plumb any cotangent bundles (of oriented manifold). 
		Second, \cite{Karabas-Lee25} allows us to use another identification map that is different from $\psi$. 
		See \cite[Section 3.4]{Karabas-Lee25}. 
		In the current article, we need the simper version of plumbing spaces defined in Definition \ref{definition plumbing space}, thus we omit the general version for convenience. 
		\item We note that a provided quiver $Q$ determines the plumbing space $P_n(Q)$, but for defining its Fukaya category, we need to specify a grading structure on it.
		Especially, if $H_1(Q)$ is nontrivial, there are multiple possible choices of a grading structure, and we need to fix one of them. 
		In \cite{Karabas-Lee25}, the choice of a grading structure can be encoded a degree function defined on the arrows of $Q$.
		More precisely, by fixing a degree for each arrows of $Q$, one can fix a grading structure of $P_n(Q)$. 
		For simplicity, in the current article, for a provided quiver $Q$, we assign a degree $0$ to every arrow of $Q$.
		We also note that by doing so, $P_n(Q)$ admits a grading structure that satisfies the {\em non-positively graded condition}, discussed in \cite{Jeong-Karabas-Lee26}.
		Thus, for the plumbing space $P_n(Q)$, we can apply the result of \cite{Jeong-Karabas-Lee26}. 
	\end{enumerate}
	 
\end{remark}

\subsubsection{Embedding of $P_n(Q)$ to $P_{n+1}(Q)$}
\label{subsubsection embbding of plumbing spaces}
Provided a quiver $Q$, we study a relation between $P_n(Q)$ and $P_{n+1}(Q)$ in Section \ref{subsubsection embbding of plumbing spaces}, which was studied in \cite[Section 7.1]{Lee24}.
Especially, we will see that $P_n(Q)$ can be embedded into $P_{n+1}(Q)$ as a symplectic submanifold. 
Or more generally, we can see that $P_n(Q)$ can be embedded into $P_{n+k}(Q)$ for all $k \geq 1$. 

We start with the simplest case, i.e., $T^*S^n$ and $T^*S^{n+1}$, which is the plumbing of the simplest quiver having one vertex and no arrows. 
We first note that 
\[S^n = \left\{(x_0, \dots, x_n) \in \mathbb{R}^{n+1}: \sum_{i=0}^n x_i^2 =1\right\} \subset S^{n+1} = \left\{(x_0, \dots, x_n, y) \in \mathbb{R}^{n+2} : \sum_{i=0}^n x_i^2 + y^2 =1\right\}.\]
Moreover, if $\eta:S^{n+1} \to S^{n+1}$ is defined by 
\[\eta(x_0, \dots, x_n, y) = (x_0, \dots, x_n, -y),\]
then $S^n$ is the set of fixed-points of $\eta$. 
On the cotangent bundle $T^*S^{n+1}$, the pull-back $\eta^*$ is an exact symplectic automorphism and $T^*S^n$ is given as the set of fixed-points of $T^*S^{n+1}$. 
Thus, $T^*S^n$ is a symplectic submanifold of $T^*S^{n+1}$. 

We recall that $P_{n+1}(Q)$ is constructed by plumbing cotangent bundles 
\[\left\{T^*S_v : v \text{  is a vertex of  } Q \right\},\] 
where $S_v \simeq S^{n+1}$. 
On each cotangent bundle $T^*S_v$, one can define an involution $\eta_v$ as we did before. 
Moreover, in the identified plumbing regions, one can choose $\eta_v$ such that $\eta_v$ commutes the identification map $\psi$ given in Equation \eqref{eqn identification map}, i.e., 
\[\eta_w \circ \psi = \psi \circ \eta_v.\]
Thus, one can construct an involution $\eta:P_{n+1}(Q) \to P_{n+1}(Q)$ such that the set of fixed-points of $\eta$ is the plumbing space of the fixed-points sets of $\eta_v$, i.e., the plumbing space $P_n(Q)$. 

Similarly, one can find an involution $\eta:S^{n+k} \to S^{n+k}$ whose set of fixed-points is $S^n$. 
And, by the same logic, we can find an involution $\eta:P_{n+k}(Q) \to P_{n+k}(Q)$ whose set of fixed-points is $P_n(Q)$, for any $n, k \in \mathbb{N}$. 
Especially, we are interested in the embedding of $P_1(Q)$ to $P_n(Q)$ for $n \geq 2$, so we will give a name for it. 
\begin{definition}
	\label{definition embedding}
	Provided a quiver $Q$, let $\eta_{(Q,n)}: P_n(Q) \to P_n(Q)$ denote the involution of $P_n(Q)$, whose set of fixed-points is $P_1(Q)$. 
	Also, let $\iota_{(Q.n)}$ denote the embedding of $P_1(Q)$ to $P_n(Q)$. 
\end{definition} 

For convenience, if there is no chance of confusion, we simply say that $P_1(Q) \subset P_n(Q)$ for $n \geq 1$, without using $\iota_{(Q,n)}$

\subsubsection{Penner type symplectic automorphisms}
\label{subsubsection Penner type symplectic automorphisms}
Now, we define our examples, symplectic automorphisms of Penner type. 
These examples are defined on plumbing spaces $P_n(Q)$ whose plumbing patterns $Q$ are {\em bipartite}.
\begin{definition}
	\label{definition bipartite}
	A finite quiver $Q$ is {\em bipartite} if $Q$ satisfies the following:
	\begin{itemize}
		\item the vertex set $V(Q)$ of $Q$ can be partitioned into two disjoint subsets 
		\[V(Q) = V_+(Q) \sqcup V_-(Q)\]
		such that every arrow has a tail in $V_+(Q)$ and a head in $V_-(Q)$;
		\item there exist at least one edge;
		\item $Q$ is connected.
	\end{itemize}
\end{definition}

\begin{remark}
	\label{rmk bipartite}
	Definition \ref{definition bipartite} is slightly stronger than the usual definition of {\em bipartite graph}. 
	The usual definition requires that every edge of a bipartite graph joins a vertex in $V_+(Q)$ and a vertex in $V_-(Q)$, and there is no direction condition. 
	In the rest of Section \ref{subsection symplectic automorphisms of Penner type}, we work with the extra direction condition for simplicity.
	For example, the direction condition simplifies the grading structure on $P_n(Q)$ (see Remark \ref{rmk a general plumbing}).
	Especially, $Q$ has no degree $0$ cycles, and induces that the zero sections $\{S_v : v \in V(Q)\}$ generate the compact Fukaya category of $P_n(Q)$ by \cite{Jeong-Karabas-Lee26}. 
\end{remark}

Especially, we note that if $Q$ is bipartite, then $Q$ cannot have a self-loop, i.e., an arrow whose head and tail are the same vertex. 
Thus, if $Q$ is bipartite, then for every vertex $v$ of $Q$, there exists an embedded Lagrangian sphere $S_v \subset P_n(Q)$. 
\begin{definition}
	\label{definition Penner type}
	Let $Q$ be a bipartite quiver whose vertex set $V(Q)$ is partitioned into $V(Q) = V_+(Q) \sqcup V_-(Q)$. 
	A symplectic mapping class $f$ of $P_n(Q)$ is of {\em Penner type} if it satisfies the following: 
	\begin{itemize}
		\item either $f$ or $f^{-1}$ has a representative, denoted as $f$ of $f^{-1}$ again, such that it is a product of positive Dehn twists along $S_v$ for $v \in V_+(Q)$ and negative Dehn twists $S_w$ along $w \in V_-(Q)$; 
		\item every $v \in V(Q)$ appears in the above product representation of $f$ or $f^{-1}$. 
	\end{itemize}
\end{definition}

From the definition, the following is an easy observation:
for a vertex $v \in V(Q)$, one can choose a Dehn twist $\tau_v$ along $S_v$ such that $\tau_v$ commutes with $\eta_{(Q,n)}$ defined in Definition \ref{definition embedding}.
Thus, if $f$ is a symplectic mapping class of $P_n(Q)$, which is of Penner type, then $f$ has a representative commuting with $\eta_{(Q,n)}$. 
It implies that the restriction of $f$ onto $\iota_{(Q,n)}\left(P_1(Q)\right) \subset P_n(Q)$ is a symplectic automorphism of the symplectic submanifold  $\iota_{(Q,n)}\left(P_1(Q)\right)$.

We note that $P_1(Q)$ is a surface, i.e., a smooth manifold of real dimension $2$. 
The mapping class group of a surface has been studied extensively. 
We can use the well-studied surface mapping class group theory to study the restriction $\overline{f}:=f|_{\iota_{(Q,n)}\left(P_1(Q)\right)}$.
Especially, from the construction of $f$, one can easily see that $\overline{f}$ is a pseudo-Anosov surface mapping class since it is constructed by {\em Penner construction} \cite{Penner88}:
note that $\overline{f}$ is also a product of positive and negative Dehn twists along circles in $P_1(Q)$, and $P_1(Q)$ is constructed by plumbing $T^*S_v^1$. 
It induces that $\overline{f}$ is constructed by Penner construction. 

Using this, one can define the following notation:
\begin{definition}
	\label{definition stretching factor}
	Let $Q$ be a bipartite quiver and $f$ be a Penner type mapping class of $P_n(Q)$. 
	As defined above, $\overline{f}$ is a pseudo-Anosov surface mapping class defined on the $2$-dimensional symplectic submanifold $\iota_{(Q,n)}(P_1(Q))$. 
	The {\em stretching factor of $f$}, denoted as $\lambda(f)$, is defined to be that of a surface mapping class $\overline{f}$. 
\end{definition}

For more details on the stretching factor of a surface mapping class, we refer the reader to \cite{Farb-Margalit12}.
And, we note that $\lambda(f) >1$ since every pseudo-Anosov surface mapping class has a positive stretching factor. 

We are ready to state the main result of Section \ref{subsection symplectic automorphisms of Penner type}.
\begin{prop}
	\label{prop Penner type}
	Let $Q$ be a bipartite quiver, and let $n \geq 3$. 
	Let $f$ be a Penner type symplectic mapping class of $P_n(Q)$. 
	We use the $\mathbb{F}_2$-coefficient for Lagrangian and fixed-point Floer theories.  
	Then, 
	\[\hfl(f) = \hcat(f) = \log \lambda(f).\]
\end{prop}

\subsubsection{Proof of Proposition \ref{prop Penner type}}
\label{subsubsection proof of Proposition}
In the proof, we use the same notation as the other parts of Section \ref{subsection symplectic automorphisms of Penner type}. 
And, since it is not a main result, but examples of the current article, we do not write all details. 

Our main strategy is to show the following series of equalities/inequalities:
\begin{gather}
	\label{eqn Penner type}
	\log \lambda(f) \stackrel{(i)}{\leq} \hfl(\overline{f}) \stackrel{(ii)}{\leq} \hfl(f) \stackrel{(iii)}{\leq} \hcat(f) \stackrel{(iv)}{\leq}  \log\lambda(f).
\end{gather} 

Let us start with (iii).
It is the simplest part since it is a direct result of Theorem \ref{thm intro}.

To prove (ii), one can apply the proof of \cite[Theorem 5.4]{Kim-Kwon26}.
We note that \cite[Theorem 5.4]{Kim-Kwon26} proved (ii) when $Q$ is the $A_k$-Dynkin diagram. 
The same proof works for a general bipartite quiver $Q$.
The key point is to show that the normal bundle of $P_j(Q)$ in $P_{j+1}(Q)$ is unitary trivial.
In the proof of \cite[Theorem 5.4]{Kim-Kwon26}, it is enough to show that $c_1(P_{j+1}(Q)) =0$. 

We note that $c_1\left(P_n(Q)\right) \in H^2\left(P_n(Q);\mathbb{Z}\right)$.
If $n >2$, we can easily see that Mayer--Vietoris sequence gives us $H^2\left(P_n(Q);\mathbb{Z}\right) = 0$, so the first Chern class vanishes. 
Similarly, if $n=1$, then $P_1(Q)$ is a Weinstein manifold of dimension $2$, and it has an one-dimensional skeleton. 
It implies that $H^2\left(P_1(Q);\mathbb{Z}\right)$ also vanishes.

For the case of $P_2(Q)$, we recall that $P_1(Q)$ can be seen as a symplectic submanifold of $P_2(Q)$ by Section \ref{subsubsection embbding of plumbing spaces}.
Conversely, one can construct $P_2(Q)$ from $P_1(Q)$ by the following procedure:
$P_1(Q)$ has a Lagrangian spheres $S_v^1$ for all $v \in V(Q)$.
Then, we consider $P_1(Q) \times \mathbb{C}$, which is a subcritical Weinstein manifold of dimension $4$. 
Along its asymptotic boundary, one can find two Legendrian spheres for each $v \in V(Q)$, which are written as 
\[S_v^1 \times \left\{\pm \infty\right\} \subset P_1(Q) \times \partial_\infty \mathbb{C} \subset \partial_\infty \left(P_1(Q) \times \mathbb{C}\right).\]
Then, $P_2(Q)$ is obtained by attaching critical Weinstein handles of real dimension $4$, $\{H_v^\pm : v \in V(Q)\}$, along $S_v^1 \times \left\{\pm \infty\right\}$, i.e., 
\begin{gather}
	\label{eqn suspension}
	P_2(Q) = \left(P_1(Q) \times \mathbb{C}\right) \cup \left(\bigcup_{v \in V(Q)} H_v^+\right) \cup \left(\bigcup_{v \in V(Q)} H_v^-\right).
\end{gather}
Then, applying Whitney sum formula for Chern classes, we have $c_1\left(P_1(Q) \times \mathbb{C}\right) = 0$.
And applying the result of \cite{Gompf98}, which explains the effect of attaching critical Weinstein/stein handles of dimension $4$, one can see that attaching Weinstein critical handles $\{H_v^\pm\}$ does not change the vanishing first Chern class. 
Thus, $c_1\left(P_2(Q)\right) = 0$. 

If $P_1(Q)$ has a genus $g \geq 2$, then (i) is a simple corollary of \cite{Felshtyn12}.
However, $P_1(Q)$ can have a single genus.
Below, we prove (i) for a general $P_1(Q)$, without assuming any genus condition. 

To prove it, we apply the results of \cite{Jiang96, Cotton-Clay09} (that studied fixed points of surface automorphisms) to $\overline{f}: P_1(Q) \to P_1(Q)$. 
Before going further, we note that $P_1(Q)$ has negative Euler characteristic if $Q$ is bipartite, so we can apply the results of \cite{Jiang96, Cotton-Clay09}. 

In \cite{Cotton-Clay09}, the author took a perturbation of $\overline{f}$, see \cite[Section 4.2]{Cotton-Clay09}.
The purpose of the perturbation is to make $\overline{f}$ to have nondegenerate periodic points and have a controllable behavior in a neighborhood of boundary components.
Let $\overline{f}_{sm}$ denote the perturbed $\overline{f}$.
Then, \cite[Theorem 1.4]{Cotton-Clay09} says that 
\[\hfl(\overline{f}) = \hfl(\overline{f}_{sm}).\]

After the perturbation, Cotton-Clay \cite{Cotton-Clay09} classified fixed points of $\overline{f}_{sm}$ into three types, see \cite[Section 4.3]{Cotton-Clay09}.
The first two types happen only if $\overline{f}_{sm}$ is a {\em reducible} mapping class, and thus in our case, we only care the third type. 

The fixed points of the third type in \cite{Cotton-Clay09} are classified again into four sub-classes.
To describe each sub-class, we note that 
\begin{enumerate}
	\item every pseudo-Anosov surface mapping class $\left[\overline{f}\right]$ (that is the same mapping class as$\left[\overline{f}_{sm}\right]$) is associated to a {\em singular} foliation $\mathcal{F}$, and every positive power $\left[\overline{f}^m\right] = \left[\left(\overline{f}^m\right)_{sm}\right]$ (that is pseudo-Anosov again) is associated to the same singular foliation $\mathcal{F}$. 
	\item our surface $P_1(Q)$ can be seen as a Weinstein domain of dimension $2$, so it is a surface with boundary components. 
\end{enumerate}

Now, we describe the four sub-classes, which we will denote by sub-classes (a)--(d) for convenience. 
Fixed points of $\overline{f}_{sm}$ in the sub-class (a) are ones {\em not} associated with any singular points of $\mathcal{F}$ or boundary components of $P_1(Q)$.
Or equivalently, we can say that the fixed points are the regular, interior fixed points of $\overline{f}$.  
The sub-classes (b) and (c) are associated with singular points of $\mathcal{F}$, especially, the type (b) is associated to an unrotated singular point while the subclass (c) is associated to a rotated singular point. 
Fixed points classified as the sub-class (d) are ones coming from the boundary components. 

Let $a(\overline{f}_{sm}), b(\overline{f}_{sm}), c(\overline{f}_{sm}), d(\overline{f}_{sm})$ denote the number of fixed points of $\overline{f}_{sm}$ in sub-classes (a)--(d), respectively. 
In \cite{Cotton-Clay09}, the author shows that $b(\overline{f}_{sm}), c(\overline{f}_{sm}), d(\overline{f}_{sm})$ are constants determined by $\mathcal{F}$. 
Then, by (2), one can easily see that for every $m \in \mathbb{N}$, 
\begin{gather}
	\label{eqn bounded}
	b\left((\overline{f}^m)_{sm}\right), \quad c\left((\overline{f}^m)_{sm}\right), \quad d\left((\overline{f}^m)_{sm}\right) \leq \text{  a constant}.
\end{gather}
Since we have a singular foliation $\mathcal{F}$, whose singular points and punctures are not depending on $m$.
As $m$ varies, a singular point could be either rotated/unrotated, thus the above numbers could vary. 
However, they are bounded.   

Also, \cite{Cotton-Clay09} showed that each fixed point in the sub-class (a) contributes an independent class of $\HF_{fixed}(\overline{f})$.
Thus, we have 
\begin{gather}
	\label{eqn lower bound}
	\hfl(\overline{f}) \geq \limsup_{m \to \infty} \frac{1}{m} \log_+ a\left((\overline{f}^m)_{sm}\right).
\end{gather}

Jiang \cite{Jiang96} proved that, for a pseudo-Anosov $\overline{f}$, the number of fixed points of $\overline{f}^m$ is asymptotic to $\lambda^m$, where $\lambda^m$ is the stretching factor of $\overline{f}$. 
More specifically, \cite{Jiang96} considered the {\em Nielsen number} of $\overline{f}$, which is defined to be the number of {\em fixed-point classes}.
We note that \cite{Cotton-Clay09} used a different terminology, {\em Nielsen classes} to indicate the fixed-point classes. 

We remark three facts. 
First, if two surface automorphisms represent the same mapping class, then they have the same Nielsen number. 
Especially, $\overline{f}^m$ and $\left(\overline{f}^m\right)_{sm}$ have the same number.
Second, since Nielsen number counts the equivalence classes, Nielsen number should be smaller than the number of fixed-points. 
It implies that the Nielsen number of $\overline{f}^m$, or equivalently that of $\left(\overline{f}^m\right)_{sm}$, is smaller than
\[a\left((\overline{f}^m)_{sm}\right)+ b\left((\overline{f}^m)_{sm}\right)+ c\left((\overline{f}^m)_{sm}\right)+ d\left((\overline{f}^m)_{sm}\right).\]
Third, \cite{Jiang96} showed that the Nielsen number of $\overline{f}^m$ is asymptotic to $\lambda^m$ if $\overline{f}$ is pseudo-Anosov and if $\lambda$ is its stretching factor. 
These three, together with inequalities in \eqref{eqn bounded}and \eqref{eqn lower bound}, imply that the inequality given in \eqref{eqn lower bound} shows that (i) holds, i.e.,  
\[\log \lambda(f) \leq \hfl(\overline{f}).\]

The final step is to show (iv). 
To prove it, we recall the computation of $\hfl(f)$.
As mentioned in Remark \ref{rmk bipartite}, since $Q$ is bipartite and since we fix a specific grading structure on $P_n(Q)$, \cite{Jeong-Karabas-Lee26} gives a specific generating Lagrangians of the compact Fukaya category.
The generating Lagrangians are $\{S_v : v \in V(Q)\}$, where $S_v$ is the zero section of $T^*S_v$ of which $P_n(Q)$ consists. 
On the other hands, \cite{Karabas-Lee25} showed that there exits a collection of generating Lagrangians $\{L_v : v \in V(Q)\}$ for the wrapped Fukaya category: $L_v$ is a cotangent fiber of $T^*S^v$.

Since $S_v$ and $L_v$ are the zero section and a cotangent fiber of $T^*S_v$, it is easy to check that $P_n(Q)$ satisfies the conditions of \cite[Theorem 4.15]{Bae-Choa-Jeong-Karabas-Lee25}.
The cited theorem provides an algorithm computing the categorical entropy of $f$.
The algorithm is explained in \cite[Section 5]{Bae-Choa-Jeong-Karabas-Lee25} in detail, and we briefly describe the algorithm below. 

The key point of the algorithm is to assign a square matrix $M_f$ to a Penner type symplectic automorphism $f$.
The size of $M_f$ is $|V(Q)|$-by-$|V(Q)|$, and for $v, w \in V(Q)$, the $(v,w)$-entry of $M_f$ is defined by
\begin{gather}
	\label{eqn matrix}
	\dim \HF^*\left(f(S_w),L_v\right).
\end{gather}
And the conclusion is that the log of {\em the maximum of absolute values of the complex eigenvalues of $M_f$ (or simply the spectral radius of $M_f$) is the categorical entropy of $f$, i.e., $\hcat(f)$}.

We note that Equation \eqref{eqn matrix} defines the matrix $M_f$, but it is not a practical way to compute $M_f$.  
There is an alternative way of computing $M_f$. 
For the alternative way, we can define a matrix $M_v$ for each $v \in V(Q)$ as follows: the $(v_1,v_2)$-component of $M_v$ is given as 
\[\dim \HF^* \left( \tau_v(S_{v_2}), L_{v_1}\right) \text{  if  } v \in V_+(Q) \text{  and  } \dim \HF^* \left( \tau_v^{-1}(S_{v_2}), L_{v_1}\right) \text{  if  } v \in V_-(Q),\]
if $\tau_v$ denotes the Dehn twist along $S_v$. 
Then, by a simple computation, we have the $(v_1,v_2)$-component of $M_v$ is 
\begin{gather}
	\label{eqn matrix2}
	\left(M_v\right)_{(v_1,v_2)} = \begin{cases}
	1 &\text{  if  } v_1 = v_2, \\
	\text{the number of arrows connecting} v_1 \text{  and  } v_2 &\text{  if  } v = v_1 \neq v_2, \\
	0 &\text{  otherwise}.
\end{cases}
\end{gather}

A  simple degree argument shows that $M_f$ is given as a product of $\{M_v : v \in V(Q)\}$.
For example, if $f = \tau_v \circ \tau_w^{-1}$, then $M_f = M_v \circ M_w$.
For more details, see \cite[Section 5]{Bae-Choa-Jeong-Karabas-Lee25}, and we will give a specific example in Section \ref{section counting conjugacy classes in symplectic mapping class groups}. 

To compare $\hcat(f)$, that is given as the log of the spectral radius of $M_f$, and the stretching factor of $\overline{f}$, we observe that 
\[H_1\left(P_1(Q);\mathbb{R}\right) = \mathbb{R} \langle S_v^1 : v \in V(Q) \rangle \oplus H_1\left(Q\right),\]
by Mayer--Vietoris sequence, where $S_v^1$ is the zero section inside $P_1(Q)$.
Then, we study the linear map $\overline{f}: P_1(Q) \to P_1(Q)$ induces on $H_1\left(P_1(Q)\right)$.

Let $\overline{f}_*$ denote the induced linear map on the first homology.
Since $\overline{f}$ can be seen as a product of positive/negative Dehn twists, we can see that $\overline{f}_*$ sends a subspace $\mathbb{R} \langle S_v^1 : v \in V(Q) \rangle$ to itself. 
Moreover, by using the basis $\{S_v^1 : v \in V(Q)\}$, we can see that over the subspace $\overline{f}_*$ acts as a matrix $M_f$, after giving a proper orientation for $S_v^1$:
it is easy to check that the Dehn twist along $S_v^1$ acts on the subspace as $M_v$ if $v \in V_+(Q)$, and the inverse of the Dehn twist along $S_v^1$ also acts as $M_v$ if $v \in V_-(Q)$. 
From this, one can show that
\begin{gather}
	\label{eqn (iv)-1}
	\hcat(f) = \log \left(\text{the spectral radius of  } M_f\right) \leq \log \left(\text{  the spectral radius of  } \overline{f}_*\right).
\end{gather} 

We note that, as stated in \cite[Example 2]{Jiang96}, there exists $g:P_1(Q) \to P_1(Q)$ such that $g$ and $\overline{f}$ are homotopic to each other and such that the topological entropy of $g$ is the log of the stretching factor of $\overline{f}$, i.e., $\log \lambda(f)$. 
We also note that if $g_*$ denotes the linear map that $g$ induces on the first homology, then $g_* = \overline{f}_*$. 
Finally, by \cite[Section 2.7(B)]{Gromov87}, we have 
\begin{gather}
	\label{eqn (iv)-2}
	\log \text{the spectral radius of  } g_* \leq \text{  the topological entropy of  } g = \log \lambda(f).
\end{gather}
Combining inequalities in \eqref{eqn (iv)-1} and \eqref{eqn (iv)-2}, we have (iv), i.e., $\hcat(f) \leq \log \lambda(f)$.

\begin{remark}
	\label{rmk conditions}
	Before ending the example subsection, we briefly explain the condition in Proposition \ref{prop Penner type}. 
	Proposition \ref{prop Penner type} uses $\mathbb{F}_2$-coefficient and assumes that $n \geq 3$. 
	Since \cite{Kim-Kwon26} uses $\mathbb{F}_2$-coefficient, we need the coefficient condition.
	Similarly, to apply \cite{Bae-Choa-Jeong-Karabas-Lee25}, we need the dimension condition, i.e., $n \geq 3$. 
\end{remark}

\section{Counting conjugacy classes in symplectic mapping class groups}
\label{section counting conjugacy classes in symplectic mapping class groups}
For the reader's convenience, we start this section by reviewing notation defined in Section \ref{section introduction}, which we needed to state Theorem \ref{thm application intro}.
We first review the notion of {\em symplectic mapping class group} and {\em mapping class group}:
for a Liouville manifold $W$,
\begin{align*}
	&\textit{Symp}(W):= \left\{ f:W \to W : f \text{  is a compactly supported, exact symplectic automorphism} \right\}, \\
	&\textit{Diff}^+(W):= \left\{f:W \to W : f \text{  is an orientation preserving diffeomorphism} \right\}.
\end{align*}
Then, there exist natural $C^\infty$-topology on the above two sets. 
And the {\em symplectic mapping class group} and the {\em mapping class group} of $W$ are defined to be 
\[\pi_0 \left(\textit{Symp}(W)\right) \text{  and  } \pi_0 \left(\textit{Diff}^+(W)\right).\]

Since every symplectic automorphism is a diffeomorphism, there exists a forgetful map from $\pi_0 \left(\textit{Symp}(W)\right)$ to $\pi_0 \left(\textit{Diff}^+(W)\right)$.
The {\em smoothly trivial symplectic mapping class group} of $W$, $K(W)$, is defined to be the kernel of the forgetful map. 

We consider the conjugate relation in $\pi_0 \left(\textit{Symp}(W)\right)$, i.e., $f \sim g$ if and only if there exists $h \in \pi_0 \left(\textit{Symp}(W)\right)$ such that $f = hgh^{-1}$. 
Using the equivalence relation, we define the following
\[C(W):= \faktor{\pi_0 \left(\textit{Symp}(W)\right)}{\sim}, \quad C\left(W;H\right) := \left\{ [f] \in C(W) : \hfl(f) <H\right\}.\]
In other words, $C(W)$ is the group of conjugate classes, and $C(W;H)$ is a subset of $C(W)$ having bounded Floer-theoretic entropy.
We note that since Floer-theoretic entropy is an invariant up to the conjugate relation, $C(W;H)$ is well-defined.
Similarly, using the fact that $K(W)$ is a normal subgroup, we define the following:
\[C_{tri}(W):= \faktor{K(W)}{\sim} , \quad  C_{tri}\left(W;H\right) := \left\{ [f] \in C_{tri}(W) : \hfl(f) < H\right\}.\] 

Now, we re-state Theorem \ref{thm application intro}.
\begin{thm}[= Theorem \ref{thm application intro}]
	\label{thm counting}
	Let $W$ be a Liouville manifold of dimension $\geq 6$ such that there exist two exact Lagrangian spheres intersecting transversally at one point. 
	\begin{enumerate}
		\item  $\left\vert C\left(W;H\right)\right\vert \geq e^H+O(1)$. 
		Especially, there exists at least exponentially many conjugacy classes in the $\pi_0 \left(\textit{Symp}(W)\right)$, having Floer-theoretic entropy $<H$. 
		\item If $\dim W = 4n$, then $\left\vert C_{tri}\left(W;H\right) \right\vert \geq \tfrac{1}{64}e^H - O(1)$. 
		Especially, there exists at least exponentially many conjugacy classes in the $K(W)$, having Floer-theoretic entropy $<H$.
	\end{enumerate}
\end{thm}

\subsection{Proof of Theorem \ref{thm counting}}
\label{subsection proof of the counting theorem}
In the current subsection, we prove Theorem \ref{thm counting}. 

To prove the first half, we recall that 
\begin{itemize}
	\item by definition, $\hfl$ is a function from $\pi_0\left(\textit{Symp}(W)\right)$ to $\mathbb{R}_{\geq 0}$;
	\item moreover, since $\hfl$ is invariant up to the conjugate relation, $\hfl$ is also a function on $C(W)$.
\end{itemize}
Thus, to prove Theorem \ref{thm counting} (1), it is enough to show that 
\begin{gather}
	\label{eqn counting 1}
	\left\vert \hfl\left(C(W)\right) \cap [0,H) \right\vert = \left\vert \hfl\left(\pi_0 \left(\textit{Symp}(W)\right)\right) \cap [0,H) \right\vert \geq e^H+O(1).
\end{gather}

Moreover, by Theorem \ref{thm local to global property}, if $V$ is a subdomain of $W$ up to completion, it is enough to prove \eqref{eqn counting 1} for $V$, instead of $W$. 
Thus, Lemma \ref{lem subdomain} makes it enough to prove \eqref{eqn counting 1} for the case of $A_2$-Milnor fiber.

\begin{lem}
	\label{lem subdomain}
	Let $(W, \lambda)$ be a Liouville manifold such that there exists two Lagrangian spheres intersecting transversally at one point. 
	Then, $W$ contains $A_2$-Milnor fiber as a subdomain up to completion, after taking a proper Liouville homotopy.
\end{lem}
\begin{proof}
	Let $L_1, L_2$ be two Lagrangian spheres of $W$ intersecting transversally at one point. 
	And for the proof, let $n$ be the dimension of $L_i$.
	
	As mentioned \cite[Section 2.1]{Smith-Wemyss23}, the union $L_1 \cup L_2$ has a small compact neighborhood $P$ such that $(P, d\lambda|_P)$ is symplectically equivalent to a Weinstein domain $(P_0,\lambda_0)$ whose completion is the plumbing of two $T^*S^n$, or equivalently, $A_2$-Milnor fiber.
	Formally, we have a symplectic embedding 
	\[j: \left(P_0, d \lambda_0\right) \hookrightarrow \left(W, d\lambda\right).\]
	
	We note that the Liouville one-form $\lambda_0$ is one defined in the $A_2$-Milnor fiber. 
	Thus, $j^*\lambda$ is not equal to $\lambda_0$.
	The goal of the proof is that by taking a proper Liouville homotopy of $\lambda$, one can make $j$ as a Liouville embedding, i.e., $j^*\lambda'=\lambda_0$ where $\lambda'$ is a modification of $\lambda$. 
	
	To do that, we first measure the difference of two one-forms $\left(j^*\lambda - \lambda_0\right).$
	Since $H^1(P_0)=0$, 
	\[j^*\lambda - \lambda_0 = dF, \text{  or equivalently  } j^*\lambda -dF = \lambda_0, \]
	for some $F: P_0 \to \mathbb{R}$. 
	
	Using $P_0 \stackrel{j}{\simeq} P \subset W$, we can find a compactly supported function 
	\[G: W \to \mathbb{R},\]
	such that $F = G \circ j$. 
	Then, we set $\lambda' = \lambda - dG$, and we have 
	\[j^* \lambda' = j^* \lambda - d \left(G \circ j\right) = \lambda_0.\]
\end{proof}

We note that the $A_2$-Milnor fiber is known to be equivalent to the plumbing space whose plumbing pattern is the $A_2$-Dynkin diagram. 
Using the notion in Section \ref{subsection symplectic automorphisms of Penner type}, let $P_n(A_2)$ denote the $A_2$-Milnor fiber. 

The advantage of seeing the $A_2$-Milnor fiber as a plumbing space is that symplectic automorphisms of Penner type are defined on the plumbing space. 
Thus, we have a large class of symplectic automorphisms, whose Floer-theoretic entropy can be computed by Proposition \ref{prop Penner type}.

For further discussion, we set the following notation: 
let $A_2$ denote the $A_2$-Dynkin quiver. 
There are two vertices, labeled $v$ and $w$, and we set $v \in V_+(A_2)$ and $w \in V_-(A_2)$. 
Thus, if $\tau_v$ and $\tau_w$ denote Dehn twists associated to $v$ and $w$, then a Penner type symplectic automorphism (or its inverse) is a product of $\tau_v$ and $\tau_w^{-1}$.  
Then, Equation \eqref{eqn matrix2} defines two matrices $M_v$ and $M_w$, associated to $\tau_v$ and $\tau_w^{-1}$, i.e., 
\[M_v = \begin{bmatrix}
	1 & 1 \\ 0 & 1
\end{bmatrix} \text{  and  } M_w = \begin{bmatrix}
	1 & 0 \\ 1 & 1
\end{bmatrix}.\]
Moreover, we could assign a matrix $M_f$ to a Penner type symplectic automorphism $f$ as follows: 
for example, if $f = \tau_v \circ \tau_w^1$, then 
\[M_f = M_v \cdot M_w.\]
As explained in Section \ref{subsection symplectic automorphisms of Penner type}, \cite{Bae-Choa-Jeong-Karabas-Lee25} proved that $\hcat(f)=\hfl(f)$ is the log of the spectral radius of $M_f$. 

We observe the following properties of $M_f$:
\begin{itemize}
	\item since $\det M_v = \det M_w =1$, $\det M_f = 1$;
	\item every components of $M_f$ is non-negative integers and $\tr \geq 3$;
	\item this implies that the characteristic polynomial of $M_f$ is 
	\[\det (M_f - t I_2) = t^2 - \tr t +1.\]
\end{itemize}
These imply that 
\[\hfl(f) = \log \left(\frac{\tr + \sqrt{\tr^2 -4}}{2}\right) = \arccosh \left(\tfrac{\tr}{2}\right), \text{  and  } \hfl(f) < H \iff \tr < 2 \cosh H.\]
	
We note that for every $T \geq 3$, if $f_T = \tau_1^{T-2} \tau_2^{-1}$, then 
\[M_{f_T} = \begin{bmatrix}
	T-1 & T-2 \\ 1 & 1
\end{bmatrix} \text{  and  } \mathrm{tr} \left(M_{f_T}\right) = T.\]
It proves Equation \eqref{eqn counting 1} since 
\[\left\vert\hfl\left(\left\{f_T: 3 \leq T < 2 \cosh H -3, T \in \mathbb{N}\right\}\right)\right\vert = \max \left\{ 0, \lceil 2 \cosh H -3 \rceil\right\} = e^H +O(1).\]

To prove Theorem \ref{thm counting} (2), we recall that, if $n$ is even, $\tau_v^8$ and $\tau_w^{-8}$ are smoothly isotopic to the identity, see \cite{Seidel14} and references therein. 
Thus, for every $T$, 
\[g_T = \tau_v^{8T} \circ \tau_w^{-8} \in K\left(P_n(A_2)\right).\]

As we did before, 
\[M_{g_T} = \begin{bmatrix}
	1+64T & 8T \\ 8 & 1
\end{bmatrix}, \text{  and  } \hfl\left(g_T\right) = \arccosh \left(1 +32T\right).\]
Finally, for every $T < \tfrac{1}{32} \left(\cosh H -1\right)$, $\hfl\left(g_T\right) < H$.
Thus, 
\[\left\vert\hfl\left(\left\{g_T: 3 \leq T < \tfrac{1}{32} \left(\cosh H -1\right), T \in \mathbb{N}\right\}\right)\right\vert = \max \left\{ 0, \lceil \tfrac{1}{32} \left(\cosh H -1\right) \rceil\right\} = \tfrac{1}{64}e^H -O(1)\]
proves Theorem \ref{thm counting} (2).

\begin{remark}
	\label{rmk local condition} 
	Theorem \ref{thm counting} has a local condition for Weinstein manifolds $W$: $W$ has two Lagrangian spheres intersecting transversally at one point. 
	We would like to assume that with a modified local condition, one can prove a modified Theorem \ref{thm counting}.
	For example, let us assume that $W$ has two Lagrangian spheres intersecting transversally at $k$-many points. 
	Or equivalently, there exists a neighborhood that is symplectomorphic to a plumbing space whose plumbing pattern is the Kronecker quiver $K_k$, i.e., $P_n(K_k)$.
	If $P_n(K_k)$ is a Liouville subdomain of $W$, then the same arguments can show the exponential abundance of the conjugacy classes in the symplectic mapping class group. 
	A different point from the local condition of Theorem \ref{thm counting} is that $P_n(K_k)$ can have a non-trivial first cohomology:
	in Lemma \ref{lem subdomain}, we used that $H^1(P_n(A_2)) = 0$. 
	Thus, for the modified local condition that two Lagrangian spheres intersecting transversally at $k$ points, we need another condition that their neighborhood is a Liouville subdomain, not just a symplectic submanifold.    
\end{remark}

\subsection{Spectrum of Floer-theoretic entropy}
\label{subsection spectrum of Floer-theoretic entropy}
The key idea of the proof of Theorem \ref{thm counting} is to show that the {\em spectrum} of Floer-theoretic entropy, i.e., 
\[\hfl\left(\pi_0 \left(\textit{Symp}(W)\right)\right) \text{  or  } \hfl\left(K(W)\right),\]
has exponentially many elements. 

Moreover, the proof also specifies a specific subset of the spectrum. 
More formally, we have 
\[\left\{\arccosh\left(\tfrac{T}{2}\right) : T \in \mathbb{N}_{\geq 3}\right\} \subset \hfl\left(\pi_0 \left(\textit{Symp}(W)\right)\right),\]
for the case of $\dim W \geq 6$, and 
\[\left\{\arccosh\left(1+32T\right) : T \in \mathbb{N}_{\geq 3}\right\} \subset \hfl\left(K(W)\right),\]
for the case of $\dim W = 4n \geq 6$.

It shows that gaps between two Floer-theoretic entropy values are {\em fine}.
Formally, we state Corollary \ref{cor gap}
\begin{cor}
	\label{cor gap}
	Under the setting of Theorem \ref{thm counting}, for any $r \in \mathbb{R}_{\geq 0}$, the distance between $r$ and the spectrum of Floer-theoretic entropy values satisfies that 
	\[\mathrm{dist} \left(r, \hfl\left(\pi_0 \left(\textit{Symp}(W)\right)\right)\right) \leq O(e^{-r}).\]
\end{cor}

We note that similar to Theorem \ref{thm counting}, Corollary \ref{cor gap} also holds with different local conditions, see Remark \ref{rmk local condition}.

\bibliographystyle{amsalpha}
\bibliography{Floer_vs_categorical}
\end{document}

%% file: first_step.pdf_tex
%% Creator: Inkscape inkscape 0.92.4, www.inkscape.org
%% PDF/EPS/PS + LaTeX output extension by Johan Engelen, 2010
%% Accompanies image file 'first_step.pdf' (pdf, eps, ps)
%%
%% To include the image in your LaTeX document, write
%%   \input{<filename>.pdf_tex}
%%  instead of
%%   \includegraphics{<filename>.pdf}
%% To scale the image, write
%%   \def\svgwidth{<desired width>}
%%   \input{<filename>.pdf_tex}
%%  instead of
%%   \includegraphics[width=<desired width>]{<filename>.pdf}
%%
%% Images with a different path to the parent latex file can
%% be accessed with the `import' package (which may need to be
%% installed) using
%%   \usepackage{import}
%% in the preamble, and then including the image with
%%   \import{<path to file>}{<filename>.pdf_tex}
%% Alternatively, one can specify
%%   \graphicspath{{<path to file>/}}
%% 
%% For more information, please see info/svg-inkscape on CTAN:
%%   http://tug.ctan.org/tex-archive/info/svg-inkscape
%%
\begingroup%
  \makeatletter%
  \providecommand\color[2][]{%
    \errmessage{(Inkscape) Color is used for the text in Inkscape, but the package 'color.sty' is not loaded}%
    \renewcommand\color[2][]{}%
  }%
  \providecommand\transparent[1]{%
    \errmessage{(Inkscape) Transparency is used (non-zero) for the text in Inkscape, but the package 'transparent.sty' is not loaded}%
    \renewcommand\transparent[1]{}%
  }%
  \providecommand\rotatebox[2]{#2}%
  \newcommand*\fsize{\dimexpr\f@size pt\relax}%
  \newcommand*\lineheight[1]{\fontsize{\fsize}{#1\fsize}\selectfont}%
  \ifx\svgwidth\undefined%
    \setlength{\unitlength}{311.81102362bp}%
    \ifx\svgscale\undefined%
      \relax%
    \else%
      \setlength{\unitlength}{\unitlength * \real{\svgscale}}%
    \fi%
  \else%
    \setlength{\unitlength}{\svgwidth}%
  \fi%
  \global\let\svgwidth\undefined%
  \global\let\svgscale\undefined%
  \makeatother%
  \begin{picture}(1,0.87272727)%
    \lineheight{1}%
    \setlength\tabcolsep{0pt}%
    \put(0,0){\includegraphics[width=\unitlength,page=1]{first_step.pdf}}%
    \put(0.85317329,0.03044794){\color[rgb]{0,0,0}\makebox(0,0)[lt]{\lineheight{1.25}\smash{\begin{tabular}[t]{l}$W^-$\end{tabular}}}}%
    \put(0.15907168,0.78725938){\color[rgb]{0,0,0}\makebox(0,0)[lt]{\lineheight{1.25}\smash{\begin{tabular}[t]{l}$W$\end{tabular}}}}%
    \put(0,0){\includegraphics[width=\unitlength,page=2]{first_step.pdf}}%
    \put(-0.09242722,0.14384079){\color[rgb]{0,0,0}\makebox(0,0)[lt]{\lineheight{1.25}\smash{\begin{tabular}[t]{l}$\text{Support of  } f$\end{tabular}}}}%
    \put(0,0){\includegraphics[width=\unitlength,page=3]{first_step.pdf}}%
    \put(0.73978051,0.71424127){\color[rgb]{0,0,0}\makebox(0,0)[lt]{\lineheight{1.25}\smash{\begin{tabular}[t]{l}$\Lambda_0$\end{tabular}}}}%
    \put(0.60317329,0.53044794){\color[rgb]{0,0,0}\makebox(0,0)[lt]{\lineheight{1.25}\smash{\begin{tabular}[t]{l}$\Delta$\end{tabular}}}}%
  \end{picture}%
\endgroup%